\documentclass{amsart}

\usepackage{epsfig, graphicx}
\usepackage{subfigure}

\usepackage{amsthm,amsmath,natbib,amsmath,amssymb,amsthm,amsfonts,amstext,amsbsy,amscd}
\RequirePackage[colorlinks=true,citecolor=blue,urlcolor=blue]{hyperref}

\RequirePackage[OT1]{fontenc}
\usepackage[applemac]{inputenc}
\usepackage{mathrsfs, amscd, amsfonts}
\usepackage{color}
\usepackage{dsfont} 
\usepackage{soul} 
\usepackage{graphicx} 
\usepackage[titletoc,title]{appendix}
\usepackage{marginnote}
\usepackage{lipsum}

\usepackage{mathtools}

\theoremstyle{plain}
\newtheorem{theorem}{Theorem}

\newtheorem{corollary}{Corollary}
\newtheorem{lemma}{Lemma}
\newtheorem{prop}{Proposition}
{\bf}{\rm}
\theoremstyle{remark}
\newtheorem{remark}{Remark}

\theoremstyle{fact}

\newcommand{\field}[1]{\mathbb{#1}}
\newcommand{\EE}{\field{E}}

\newcommand{\GG}{\field{G}}

\newcommand{\NN}{\field{N}}
\newcommand{\PP}{\field{P}}

\newcommand{\RR}{\field{R}}
\newcommand{\TT}{\field{T}}

\newcommand{\Ff}{{\mathcal F}}

\newcommand{\Ll}{{\mathcal L}}

\newcommand{\Nn}{{\mathcal N}}

\newcommand{\Tt}{{\mathcal T}}

\newcommand{\indicatrice}{\mathds{1}}

\newcommand{\nlim} {\underset{n \rightarrow +\infty}{{\xrightarrow{\hspace*{0.5cm}}}} }
\newcommand{\plim} {\underset{p \rightarrow +\infty}{{\xrightarrow{\hspace*{0.5cm}}}} }

\def \ep {\varepsilon}

   \newcommand\smallO[1]{
      \mathchoice
         {
            \ensuremath{\mathop{}\mathopen{}{\scriptstyle\mathcal{O}}\mathopen{}\left(#1\right)}
         }
         {
            \ensuremath{\mathop{}\mathopen{}{\scriptstyle\mathcal{O}}\mathopen{}\left(#1\right)}
         }
         {
            \ensuremath{\mathop{}\mathopen{}{\scriptscriptstyle\mathcal{O}}\mathopen{}\left(#1\right)}
         }
         {
            \ensuremath{\mathop{}\mathopen{}{o}\mathopen{}\left(#1\right)}
         }
   }

\DeclareMathAlphabet{\mathpzc}{OT1}{pzc}{m}{it}

\begin{document}
\title[Large deviations and large local densities]{A phase transition for large  values of \\
bifurcating autoregressive models}

\author{Vincent Bansaye and S. Val\`ere Bitseki Penda}

\address{Vincent Bansaye, CMAP,
 \'Ecole Polytechnique, Route de Saclay, 91128 Palaiseau, France.}

\email{vincent.bansaye@polytechnique.edu}

\address{S. Val\`ere Bitseki Penda, Institut de Math\'ematiques de Bourgogne, UMR 5584, CNRS,
Universit\'e de Bourgogne Franche- Comt\'e, F-21000 Dijon, France.}

\email{Simeon-Valere.Bitseki-Penda@u-bourgogne.fr}

\begin{abstract}
We describe the  asymptotic behavior of the number  $Z_n[a_n,\infty)$ 
of individuals with a large value  in a stable bifurcating autoregressive process.
The study of the associated first moment  $\EE(Z_n[a_n,\infty))$ 
is equivalent to the annealed  large deviation problem  $\PP(Y_n\geq a_n)$,  where 
$Y$ is an autoregressive process in a random environment and $a_n\rightarrow \infty$.\\
The population with large values  and the trajectorial  behavior   of  $Z_n[a_n,\infty)$ 
 is obtained from
 the ancestral paths associated to the large deviations of $Y$ together with its environment.
The study  of  large deviations of autoregressive processes in random environment is of independent interest and achieved first  in this paper. 
The proofs of trajectorial estimates  for  bifurcating autoregressive process  involves then a  law of large numbers for  non-homogenous trees.\\
Two regimes appear in the stable case, depending on the fact that one of the autoregressive parameter is greater than one or not.
It yields two different  asymptotic behaviors for the large local densities and maximal value of the bifurcating autoregressive process.
\end{abstract}
\maketitle



\section{Introduction}
The bifurcating autoregressive (BAR) process $X=(X_n)_{n\geq 1}$ is a model for affine random transmission of a real value  along a binary tree. To define this process, we consider a real value
 random variable $X_1$, independent of  a sequence of i.i.d bivariate random variables $\big((\eta_{2k}, \eta_{2k+1}), k\geq 1\big)$  
 with law $\mathcal{N}_{2}(0,\Gamma)$, where 
 \begin{equation*} \Gamma=\begin{pmatrix} 1 & \rho \\ \rho & 1 \end{pmatrix},\qquad \rho\in(-1,1). \end{equation*}  Then $X$ is defined inductively for $k\geq 1$ by
\begin{equation}\label{bar11} \begin{array}{ll} \quad \left\{\begin{array}{ll} X_{2k}=\alpha X_{k} + \eta_{2k} \\ \\ X_{2k+1} = \beta X_{k} + \eta_{2k+1}, \end{array} \right. \end{array}\end{equation}
where $\alpha,\beta$ are non-negative real numbers. Informally, the value $X_k$ of individual $k$ is randomly transmitted to its two offsprings $2k$ and $2k+1$ following an autoregressive process.\\
 
This model  has been introduced in the symmetric case $\alpha=\beta$  by Cowan \cite{Cowan84} and Cowan and Staudte \cite{CS86} to analyze  cell lineages data.
 It allows to study the evolution and the transmission of a trait after division, in particular it size or it growth rate. In these works, Cowan and Staudte have  focused on the study of the bacteria Escherichia Coli. E. Coli is a rod-shaped bacterium which reproduces by dividing in the middle, producing two cells. Several extensions of their model have been proposed and  we refer e.g. to the works of Basawa and Higgins \cite{BH99,BH00} and Basawa and Zhou \cite{BZ04,BZ05}, where the model of Cowan and Staudte is studied for long memory and non-Gaussian noise.    \\

 In 2005, Steward et \textit{al.} \cite{SteMadPauTad} have designed an experimental protocol which brings evidence of aging and asymetry in E. Coli. In order to study  the dynamic of evolution of values in a population of cells, taking into account  such possible asymmetry, Guyon and \textit{al.} \cite{Gu&Al} have considered and used the model of Cowan and Staudte with $(\alpha, \beta) \in (0,1)^{2}$, which means that  $\alpha$ and 
 $\beta$ may be different.
Since then, this model has been extended to more complex settings and studied from a statistical and probabilistic point of view. Bercu et al.~\cite{BdSGP09} consider an extension of BAR model with non- Gaussian noise and long memory. They use martingale approach in order to study the asymptotic behavior of least- squares estimators of unknown parameters of their model. In the following of this work, Bitseki and Djellout \cite{BD14} and Bitseki et al.~\cite{BDG11} have studied the deviation inequalities and the moderate deviations principle of unknown parameters in the BAR models introduced by Bercu et al.~\cite{BdSGP09} and Guyon \cite{G07}. Several extension of this model were proposed where the missing data are taking into account in the study of cell lineage data, see  de Saporta et al.~\cite{dSGPM11,dSGPM12,dSGPM14} and  Delmas and Marsalle \cite{DM10}. Furthermore  Bercu and Blandin \cite{BB13} and de Saporta et al.~\cite{dSGPM13}, have considered
random coefficients, while Bitseki Penda and Olivier \cite{BO16} have worked on nonlinear  autoregressive models.\\

More generally, polarisation of cells is a fundamental process in cell biology. Asymmetry at cell division have been observed and studied in different contexts, see e.g. \cite{Kimmel} for plasmids, \cite{Sinclair} for extra chromosomal DNA, \cite{K} for mitochondrias 
 and \cite{bansaye} for parasites. Asymmetry at division is  a key feature of aging, cell variability and differentiation \cite{chartier}. \\
 The BAR models turns out  to be  an interesting toy  model to explore mathematically the effect of randomness and asymmetry in the long time behavior of characteristics of cells in division processes.  In this paper, we focus on large values among the cell population. We can, describe here  how lineages and stochasticity from the Gaussian additive term interplay to explain this part of the distribution of traits when time goes to infinity.

\subsection{Stability and random cell lineage.}
To analyze the long time behavior of the bifurcating autoregressive model, the genealogy of the population is involved. It is given here by the binary tree.
Each vertex is seen as a positive integer different from 0 and it denotes one individual of the population. The initial individual is thus denoted by $1$ and for $n\in\mathbb{N}$,
\begin{equation*} \mathbb{G}_{n}=\{2^{n},2^{n}+1,\cdots,2^{n+1}-1\} \quad \text{and} \quad \mathbb{T}_{n}=\bigcup_{m=0}^{n}\mathbb{G}_{m} \end{equation*}
denote respectively the $n$-th generation and the first $(n+1)$ generations of the population. 
The collection of values $(X_i : i\in \mathbb G_n)$ of the individuals in generation $n$ is represented by the random punctual measure 
$$Z_n=\sum_{i\in \mathbb G_n} \delta_{X_i}.$$
The process $Z$ is a multitype branching process where types are real valued. The first moment of  measure  $Z_n$ satisfies the following simple many-to-one formula:
\begin{equation}
\label{MTOf}
\mathbb E\left(Z_n([a,b])\right)=2^n \mathbb P( Y_n\in [a,b]),
\end{equation}
for any real numbers $a\leq b$,
where  $(Y_n)_n$ is the autoregressive process in random environment corresponding to a uniform random path in the  binary tree. More precisely, $Y$ is defined by
\begin{equation}\label{eq:random-ar}
Y_{0} = X_1 \quad \text{and} \quad \forall n\geq 1, \quad Y_{n} =
\theta_{n} Y_{n-1} + \ep_{n},
\end{equation}
where $\Theta=(\theta_{n},n\geq 1)$ is a sequence of
i.i.d. random variables such that
$$\mathbb P(\theta_1=\alpha)=\mathbb P(\theta_1=\beta)=1/2$$
and  $(\ep_{n},n\geq 1)$ is  a sequence of
i.i.d. centered standard Gaussian random variables independent of
$(\theta_{n},n\geq 1)$. We observe that for any $n\geq 1$,
\begin{equation}\label{eq:yn}
Y_{n} = \left(\prod_{k=1}^{n} \theta_{k}\right)Y_{0} +\overline{Y}_n,  \text{ where }  \quad \overline{Y}_n= \sum\limits_{k=1}^{n} \left(\prod_{\ell=k+1}^{n} \theta_{l}\right)\ep_{k},
\end{equation}
using  the convention $\Pi_{\ell \in \varnothing}=1$.
Conditionally on $(\theta_{n},n\geq 1)$, $\overline{Y}_n$ is
a Gaussian random variable with  standard deviation
\begin{equation}\label{eq:An}
A_n=A(\theta_1,...,\theta_n)=\sqrt{\sum_{k=1}^{n} \prod_{\ell=k+1}^{n} \theta_{\ell}^{2}}.
\end{equation}
Moreover,  for each $n$, $A_n$ is distributed as $A_n^*=A(\theta_n,...,\theta_1)$ and by monotonicity the latter converges a.s. as $n\rightarrow \infty$ 
to  
$$A^*_{\infty}=\sqrt{\sum_{k=1}^{\infty} \prod_{\ell=1}^{k-1} \theta_{\ell}^{2}} \in (0,\infty].$$
If $\alpha\beta \geq1$ then $A^*_{\infty}=\infty$ a.s. and the mean  proportion of individuals whose value is in a compact set goes to $0$ as $n$ goes to infinity, i.e. $Z_n([a,b])/2^n \rightarrow 0$ a.s.
We focus here on the stable case  :
$$\alpha\beta < 1.$$
Then $(Y_n)_n$ converges in law to its unique stationary distribution $\pi$ on $\mathbb R$, which is a  mixed centered gaussian random variable with random standard deviation  distributed as $A^*_{\infty}$.
Guyon \cite{G07}  has shown in  that case that a strong law of large numbers holds : 
$$\frac{1}{2^n}Z_n([a,b])\stackrel{n\rightarrow \infty}{\longrightarrow} \pi([a,b]) \qquad \text{a.s.}$$

\subsection{Large deviations and local densities in the stable case : main results.}
We are interested here in the number of individuals with larges values  and the extremal values among the population. More precisely we study
$$Z_n([a_n,\infty))=\#\{ i \in \mathbb G_n : X_i \geq a_n\}$$ 
when $a_n\rightarrow \infty$. In particular we expect a  law of large number (or concentration)  effect, informally  
$$Z_n([a_n,\infty)) \approx \mathbb E(Z_n([a_n,\infty)))= 2^n \mathbb P(Y_n\in [a_n, \infty))$$
at the logarithm scale. \\
Proving this law of large number effect and determining the asymptotic behavior of $\mathbb E(Z_n([a_n,\infty))$
leads us to study  the  large deviation event $\{Y_n\geq a_n\}$. More precisely,  we describe  the past trajectory $(Y _i  : i\leq n)$  conditionally on this event $\{Y_n\geq a_n\}$  together with the associated environment $(\theta_i :i\leq n)$. It  allows us both to estimate the mean behavior $\PP(Y_n\geq a_n)$ and to control the subtree associated to the local density $Z_n([a_n,\infty))$ and thus correlations of the values between individuals.  This large deviation issue on the autoregressive process $Y$ is  both crucial here
for the analysis of $Z_n$ and    of independent interest. 
In the case $\alpha=\beta<1$,   $Y_n$ is a gaussian random variable and $\mathbb P(Y_n\geq a_n)$ is derived analytically. In the general stable case $\alpha \beta <1$, 
two regimes appear, which correspond to two different behaviors of the past trajectory and associated environments $\theta$. 
We obtain  the following classification
\begin{itemize}
\item[$i)$] Case $\alpha=\max\{ \alpha, \beta\} < 1$. Then, the deviation event $\{Y_n\geq a_n\}$ is achieved by selecting the most favorable environment $\alpha$ in the last generations. The deviation event  $\{Y_n\geq a_n\}$ relies on  extreme events in these last generations, which are quantified by the tail of  Gaussian random variables with standard deviation $\alpha$. We obtain that  for any $x>0$,
\begin{equation*}
\lim_{n\rightarrow+\infty}\frac{1}{xn}\log\PP\left(Y_{n}\geq
\sqrt{xn}\right) = -\frac{1-\alpha^{2}}{2}.
\end{equation*}
Moreover we prove that conditionally on $\{Y_n\geq a_n\}$, the past trajectory $(Y _{n-i}  : 0\leq i\leq n)$  is approximated   by a geometric progression
$(a_n \alpha^i :  0 \leq i\leq n)$. 
\item[$ii)$] Case $\beta< 1 <\alpha$  and $\alpha\beta<1$.  Then there exists a unique $\kappa>0$  such that
$$\alpha^{\kappa}+\beta^{\kappa}=2$$
and  now the deviation event $\{Y_n \geq a_n\}$ comes from the deviation of the environment $\{A_n\geq a_n\}$.
To estimate the probability of that latter, we  use Kesten theorem, which ensures that $\mathbb P(A_{\infty} \geq a_n) \sim ca_n^{\kappa}$ as $n\rightarrow \infty$.
We thus restrict the study here to the case when $\{A_n^* \geq a_n\}$ is comparable to $\{A_{\infty}^*\geq a_n\}$, which will correspond to the fact that  the supremum over all time of the associated random walk drifting to  $-\infty$ is reached before time $n$  when it is larger than $a_n$. That leads us to introduce $\gamma\in (0,1)$  defined by 
$$ 
\gamma :=\inf_{s\geq 0}( \alpha^{s}+\beta^{s})/2.
$$ 
Then for a sequence $a_n$ tending to infinity not too fast, we prove that $\mathbb P(Y_n\geq a_n)$ is equivalent to $a_n^{-\kappa}$. More precisely we need that  $a_n^{\kappa}$ is negligible compared to   $n^{3/2}(1/\gamma)^{n}$ as $n$ tends to infinity. In particular,  for any $\rho \geq 1$ such that $\rho^{\kappa} \in (0,1/\gamma)$, we get
\begin{equation*}
\lim_{n\rightarrow+\infty} \frac{1}{n\log(\rho)}\log\PP\left(Y_{n} \geq \rho^n \right) =-\kappa.
\end{equation*}
\end{itemize}
We expect that the regimes $\max(\alpha,\beta)=1$ and $a_n=\exp(n\rho)$ with $\rho>1/\gamma$ could be studied with an approach similar to $ii)$. As far as we seen, $\{A_n^*\geq a_n\}$ will not be comparable to $\{A_{\infty}^*\geq a_n\}$ if $\rho>1/\gamma$ and describing the behavior of $\{Y_n\geq a_n\}$  should require more work. We also mention that non gaussian noise could lead to different behavior, in particular if the  tail of $\varepsilon$ is  at the same scale as the deviation probability coming from  the environment.\\
$\newline$

We have thus obtained both the asymptotic behavior of $\mathbb E(Z_n([a_n,\infty)))= 2^n \mathbb P(Y_n\in [a_n, \infty))$ and the ancestral path of $(Y, \Theta)$ on the event   $Y_n\geq a_n$. Roughly speaking, a trajectorial version of the 
  many-to-one formula \eqref{MTOf}   allows  us to derive
 the genealogy and process  contributing to the local density 
 $Z_n([a_n,\infty))$ from the  ancestral path of $(Y, \Theta)$ before time $n$.  We obtain the following long time estimates
 \begin{itemize}
\item[$i)$] Case $\alpha=\max\{ \alpha, \beta\} < 1$. Then, for any $x\in (0,2\log(2)/(1-\alpha^{2}))$,
$$ \lim_{n\rightarrow+\infty} \frac{1}{n} \log
Z_{n}([x\sqrt{n};+\infty)) = \log(2) - x(1-\alpha^{2})/2 \qquad \text{a.s}.$$
In that case, the proof follows ideas of Biggins  for branching random walks and we consider a subpopulation of individuals which realizes the deviation and thus the large values $[x\sqrt{n},\infty)$. It is  described by a  non-homogeneous Galton-Watson process, inherited from the trajectory
of  values in last generation following the geometric paths obtained above for the deviation $\{Y_n \geq a_n\}$. The main difference with branching random walk is this time non-homogeneity since the trajectory associated with the deviation of the value (or value) is not straightline but geometric and limited to the last generations.
\item[$ii)$]  Case $\beta < 1 <\alpha$ and  $\alpha\beta<1$.
For all $\rho \geq 1$ such that $\rho^{\kappa} < \min(2,1/\gamma)$, we have
\begin{equation*}
\lim_{n\rightarrow +\infty} \frac{1}{n} \log (Z_n([\rho^n,\infty)) = \log(2) - \kappa\log (\rho) \quad \text{in probability}.
\end{equation*}
The subtree associated with the individuals in $[\rho^n,\infty)$ is very different from the previous case. Indeed, only some particular paths of the binary genealogical tree are involved in local densities, corresponding to the deviation event $\{A_n^*\geq \rho^ n\}$.  The evolution of the chain on these  paths does not deviate from its usual Gaussian autoregressive  expected  behavior.   The proof is here more involved since a  study  of the process restricted to the subtree of $\TT_n$ corresponding to the environment $\Theta_n=(\theta_i :i \leq n)$ associated to the large deviation of $Y_n$.
This phenomenon  is linked to the random environment structure inherited from the tree and  can be compared to the results in the weakly subcritical case
in  Kimmel's branching model \cite{bansaye}. In that latter,  an analogous phase transition occurs but only the mean behavior had been obtained.  The technique developed in this paper should  help  to obtain a  convergence in probability. The issue of the a.s. convergence remains open. We could expect to prove it by estimating for instance the speed of convergence in probability.
\end{itemize}

In words, two sources of randomness are combined here : the gaussian additive term of the autoregressive process and the choice of paths in the binary genealogy. In the first regime, the local densities is inherited from both randomness, in the very last generations.  The picture is quite close to branching random walks, but the fact that deviation creating local densities occurs only at the end, not starting from the beginning.
In the second regime,  only the path in the tree is involved  for large traits and the deviation starts much earlier to create these local densities. The picture is very different from branching random walks since a deterministic subtree, inherited from functional $A_n$ support the local density. The study of the autoregressive model $Y$ together with its environment $\theta$  yields a natural way to analyse this issue.

\subsection{Extremal value.}
The  maximal value in generation $n$ 
$$M_n=\max\{ X_i : i\in \mathbb G_n\}$$
has different behavior and speed. If $\alpha=\max\{ \alpha, \beta\} < 1$,
then $$ \lim_{n\rightarrow+\infty} \frac{1}{\sqrt{n}} M_n=\sqrt{\frac{2 \log(2)}{1-\alpha^{2}} }\qquad \text{a.s}.$$
If  $\beta < 1 <\alpha$, then
 $$ \lim_{n\rightarrow+\infty} \frac{1}{n} \log(M_n)=\frac{\log(2)}{\kappa} \qquad \text{in probability}.$$ 
The issue of maximal value has attracted lots of attention  and has been extensively studied for branching random walks from the works of Hammersley \cite{Hammersley74}, Kingman \cite{Kingman75}, Biggins \cite{B76} and Bramson \cite{Bramson78}. Fine  estimates have been obtained recently, which have shed light on the branching structure
of the process, see in particular Hu and Shi \cite{HS09}, Faraud et \textit{al.} \cite{FHS12}, Addario-Berry and Reed \cite{AR09} and A\"{\i}d\'ekon \cite{Aidekon13}.\\
The behavior of the maximal value is different in this model.
Technics developed for branching random walk  are used and adapted here, in particular the use of many to one formula to exploit the description of the path leading to extremal 
particules \cite{Rob} and the construction of an imbedded branching process. 
The main novelty lies in the weak regime $\beta<1< \alpha$ where a law of large numbers of a well chosen subtree of the binary tree is involved.

\subsection{Notation.} We will use the following notations. Let $(x_{n})$ and $(y_{n})$ be two sequences of real numbers. 
\begin{itemize}
\item We recall that $x_{n} = \smallO{1}$ means that $x_{n} \rightarrow 0$ as $n \rightarrow +\infty.$
\item We write $x_{n} = \smallO{y_{n}}$ if $x_{n}/y_{n} = \smallO{1}.$
\item We write $x_{n} \sim y_{n}$ if $x_{n}/y_{n} \rightarrow 1$ as $n \rightarrow +\infty$.
\item We write $x_{n} \lesssim y_{n}$ if there is a positive constant $c$, independent of $n$, such that $x_{n} \leq c y_{n}$.
\item We write $x_{n} \asymp y_{n}$ if there are two positives constants $c_{1}$ and $c_{2}$ such that $c_{1}y_{n} \leq x_{n} \leq c_{2}y_{n}$ as $n \rightarrow +\infty$.
\end{itemize}  


\section{Large deviations and ancestral paths  in the stable case}

This section has two goals. The first is to evaluate the probabilities of large deviations of random coefficients autoregressive process $(Y_{n}, n\in\NN)$. The second is to understand the asymptotic behavior of the paths of this process conditionally to the large deviations events. For simplicity and owing to the motivation of this paper (see the next section), 
we focus on the case where the random coefficients  take only two values with identical probability. Then, let $(\theta_{n},n\geq 1)$ be a sequence of i.i.d. random variables which take their values in $\{\alpha,\beta\}$ with $0 < \beta \leq \alpha$, and $\PP(\theta_{1} = \alpha) = 1/2$. Let $(\ep_{n},n\geq 1)$ be a sequence of i.i.d. centered standard Gaussian random variables independent of $(\theta_{n},n\geq 1)$. We assume without loss of generality that $Y_{0} = 0$ and  consider the process $(Y_{n},n\in \NN)$ defined by
\begin{equation*}\label{eq:random-ar2}
Y_{0} = 0 \quad \text{and} \quad \forall n\geq 1, \quad Y_{n} = \theta_{n} Y_{n-1} + \ep_{n},
\end{equation*}
in such a way that from \eqref{eq:yn}, we have
$$
Y_{n} = \sum\limits_{k=1}^{n} \left(\prod_{\ell=k+1}^{n} \theta_{\ell}\right)\ep_{k}.
$$
  
\subsection{Large deviations for strictly stable case}$\,$
We consider here the case $$\beta\leq \alpha < 1$$
and show that the (annealed) large deviation behavior of  $(Y_{n},n\in\NN)$ is given by a geometric growth in envrionment $\alpha$  in the last generations.
\begin{theorem}\label{thm:large_deviation}
Let $(a_{n})_{n\in \NN}$ be a sequence of non-decreasing real numbers which tends to infinity.
We have
\begin{equation}\label{eq:ldyn}
\lim_{n\rightarrow+\infty}\frac{1}{a_{n}^{2}}\log\PP\left(Y_{n}\geq a_{n}\right) = -\frac{1-\alpha^{2}}{2}.
\end{equation}
Moreover, for all sequence $(\ell_{n})_n$ of  integers such that $n-\ell_n\rightarrow \infty$ and $\ell_{n} = \smallO{\log{a_{n}}}$, 
\begin{equation}\label{eq:behav-theta}
\lim_{n\rightarrow +\infty} \PP\left(\sup_{\ell\in \{ 0,\cdots,\ell_{n}\}}  \left|\frac{Y_{n-\ell}}{a_{n}} - \alpha^{\ell}\right| \leq \varepsilon, \quad   (\theta_{n-\ell_n}, \ldots, \theta_n)= (\alpha, \ldots, \alpha) \,  \bigg\vert  \, Y_{n} \geq a_{n}\right) = 1
\end{equation}
for any  $\varepsilon>0$.
\end{theorem}
\begin{remark}
More generally,  one can see the proof to check that admissible sequences $(\ell_n)_n$ just need to satisfy $n-\ell_n\rightarrow \infty$
and  $\lim_{n\rightarrow \infty} 2^{\ell_{n}} \exp\left(-c (a_{n}\alpha^{\ell_{n}})^{2}\right) = 0$ for some finite constant $c$.
\end{remark}
\begin{proof}
 Let $\Ff_{n} = \sigma(\theta_{1},\cdots,\theta_{n})$ be the filtration defined by $\theta_{1},\cdots, \theta_{n}$. We recall that $\Ll(Y_{n}|\Ff_{n}) = \Nn(0,A_{n}^{2})$, the centered normal law with variance $A_{n}^{2}$, where $A_{n}$ is defined in \eqref{eq:An}. We observe that $$A_n^2\leq \sum_{k=0}^{\infty} \alpha^{2k}=\frac{1}{1-\alpha^2}$$ since $\beta\leq \alpha$. First, we have
\begin{align}
\PP\left(Y_{n} \geq a_{n}\right) &
= \EE\left[\PP\left(\Nn(0,A_{n}^{2}) \geq a_{n} \big| \Ff_{n}\right) \right] \nonumber\\ &
\leq \PP\left(\Nn\left(0, \frac{1}{1-\alpha^{2}}\right) \geq a_{n}\right)\lesssim \frac{1}{a_{n}} \exp\left(-\frac{(1-\alpha^{2})a_{n}^{2}}{2}\right),& \label{eq:upyn}
\end{align}
where the last inequality follows from the  classical upper tail inequality for standard normal distribution,  see for instance \cite{Birnbaum42}. \\
Next, we introduce a sequence  $(\ell_{n})_{n\in\NN}$ which tends to infinity such that $\ell_n\leq n $ and $\ell_{n} = \smallO{a_{n}^{2}}$. We set 
$$B_{\ell_{n}} = \{\theta_{n} = \alpha,\cdots,\theta_{n-\ell_{n}} = \alpha\}, \qquad u_{\ell_{n}}^2=\sum\limits_{k=0}^{\ell_{n}} \alpha^{2k}=\frac{1 - \alpha^{2\ell_{n} + 2}}{1-\alpha^2}$$
and using again 
 $\Ll(Y_{n}|\Ff_{n}) = \Nn(0,A_{n}^{2})$,  we have
\begin{align} 
\PP\left(Y_{n} \geq a_{n}, B_{{n}}\right) &= \frac{1}{2^{\ell_{n}+1}} \times \EE\left[\PP\left(\Nn\left(0, \alpha^{2\ell_{n} + 2} A(\theta_{1},\cdots,\theta_{n-\ell_{n}-1})^{2} + u_{\ell}^2 \right) \geq a_{n} \big| \Ff_{n-\ell_{n}-1}\right)\right]& \nonumber\\ &\geq \frac{1}{2^{\ell_{n}+1}} \times \PP\left(\Nn\left(0,u_n^2\right) \geq a_{n}\right) \nonumber\\ &\gtrsim \frac{1}{2^{\ell_{n}}} \times \frac{u_{\ell_{n}}}{a_{n}} \exp\left(- \frac{a_{n}^{2}}{2u_{\ell}^2}\right),& \label{eq:loyn}
\end{align} 
where the last inequality follows from the classical lower bounds for the tail of standard normal distribution (see for e.g. \cite{Birnbaum42}). \\
Now, from $\eqref{eq:upyn}$ and $\eqref{eq:loyn}$ we get
$$
\frac{1}{2^{\ell_{n}}}  \frac{u_{\ell_{n}}}{a_{n}} \exp\left(- \frac{a_{n}^{2}}{2u_{\ell_{n}}^2}\right) \lesssim \PP\left(Y_{n} \geq a_{n}\right)  \lesssim \frac{1}{a_{n}} \exp\left(-\frac{(1-\alpha^{2})a_{n}^{2}}{2}\right). 
$$
Finally, letting $n$ go to  infinity in the previous inequalities leads to 
\begin{equation*}
-\frac{1-\alpha^{2}}{2} \leq \liminf_{n\rightarrow+\infty} \frac{1}{a_{n}^{2}} \log\PP\left(Y_{n}\geq a_{n}\right) \leq \limsup_{n\rightarrow+\infty} \frac{1}{a_{n}^{2}}\log\PP\left(Y_{n}\geq a_{n}\right) \leq -\frac{1-\alpha^{2}}{2}
\end{equation*}
since  $\ell_{n} = \smallO{a_{n}^{2}}$ and $u_{\ell_{n}}^2\rightarrow 1/(1-\alpha^2)$, which  ends the proof of \eqref{eq:ldyn} . \\

For the second part,  we focus on the case $\beta<\alpha$, while the case   $\alpha=\beta$ is simpler. 
We write 
$$
\tau_n =\sup\{  i=1,\ldots, n :\theta_i= \beta\}
$$
with convention $\sup\varnothing=0$. We first observe that
$$
\PP\left(Y_{n} \geq a_{n}; \tau_n=n-i\right)=\PP\left(\theta_{n} = \alpha,\cdots,\theta_{n-i+1} = \alpha, \theta_{n-i} = \beta, Y_{n} \geq a_{n}\right)
$$
and
$$\mathcal L (Y_n \vert \theta_{n} = \alpha,\cdots,\theta_{n-i+1} = \alpha, \theta_{n-i} = \beta)=\mathcal L\left(\alpha^{i}\beta Y_{n-i-1} + \sum\limits_{k=0}^{i} \alpha^{k} \ep_{n-k}\right).$$
We get
\begin{eqnarray*}
&&\PP\left(Y_{n} \geq a_{n}; \tau_n=n-i\right)\\
&&\qquad =\left(\frac{1}{2}\right)^{i+1} \times  \PP\left(\alpha^{i}\beta Y_{n-i-1} + \sum\limits_{k=0}^{i} \alpha^{k} \ep_{n-k} \geq a_{n}  \right)\\
&&\qquad =\left(\frac{1}{2}\right)^{i+1} \times \mathbb E\left( \PP\left(\Nn\left(0,\frac{1-(\alpha^{2})^{i+1}}{1 - \alpha^{2}} + \left(\alpha^{i}\beta\right)^{2} A_{n-i-1}^{2}\right)\geq a_{n} \big| 
\Ff_{n-i}\right) \right)\\
&&\qquad \leq \left(\frac{1}{2}\right)^{i+1} \times \mathbb  \PP\left(\Nn\left(0,w_i^2\right)\geq a_{n} \right),
\end{eqnarray*}
where using again $A_{n-i-1}\leq 1/(1-\alpha^2)$,  $w_i$ is a non-negative real number defined by $$w_i^2=\frac{1-(\alpha^{2})^{i+1}}{1 - \alpha^{2}} + \frac{\left(\alpha^{i}\beta\right)^{2}}{1-\alpha^2}=\frac{1-\alpha^{2i}(\alpha^2-\beta^2)}{1-\alpha^2}.$$
Then
$$
\PP\left(Y_{n} \geq a_{n}; \tau_n=n-i\right)
\lesssim \left(\frac{1}{2}\right)^{i+1}\frac{1}{a_n}\exp\left(-\frac{a_n^2}{2w_i^2}\right)
$$
and summing over $i$, we get by monotonicity of $w_i$,
$$\PP\left(Y_{n} \geq a_{n},B_n^c\right)\lesssim  \sum_{i=0}^{\ell_n} \left(\frac{1}{2}\right)^{i+1}\frac{1}{a_n}\exp\left(-\frac{a_n^2}{2w_i^2}\right)\lesssim \frac{1}{a_n} e^{-a_n^2/w_{\ell_n}^2}.$$
Let $s \in \NN$ such that $\alpha^{2s} < \alpha^{2} - \beta^{2}$. Using now $\PP\left(Y_{n} \geq a_{n}\right)\geq \PP\left(Y_{n} \geq a_{n}, B_{\ell_{n} + s}\right)$ and \eqref{eq:loyn}, we get
$$
\limsup_{n\rightarrow \infty} \frac{\PP\left(Y_{n} \geq a_{n},B_n^c\right)}{\PP\left(Y_{n} \geq a_{n}\right)}\lesssim \limsup_{n\rightarrow \infty} 2^{\ell_{n} + s} e^{-a_n^2(1/w_{\ell_n}^2-1/u_{\ell_{n} + s}^2)}=0,
$$
since    $1/w_{\ell_n}^2-1/u_{\ell_{n} + s}^2\sim c\alpha^{2\ell_n}$ with $c>0$ and $w_{\ell_n}/u_{\ell_{n} + s}\rightarrow 1$ and $\ell_n= \smallO{\log{a_{n}}} $. It ensures that
\begin{equation}
\label{partieenv}
\lim_{n\rightarrow +\infty} \PP\left(B_{\ell_{n}} | Y_{n} \geq a_{n}\right) = 1.
\end{equation}
This proves a part of   \eqref{eq:behav-theta}.\\
Besides, recalling that $n-\ell_n\rightarrow \infty$,  we consider $k_n\leq n$ such that $k_n-\ell_n\rightarrow \infty$  and $k_n=\smallO{\log{a_{n}}}$. Then, \eqref{partieenv} ensures that
\begin{multline*}
\limsup_{n\rightarrow \infty} 
\PP\left(\sup_{\ell = 0,\cdots,\ell_{n}} \left|\frac{Y_{n-\ell}}{a_{n}} - \alpha^{\ell}\right| > \varepsilon \Big| Y_{n}\geq a_{n}\right) 
\leq  \limsup_{n\rightarrow \infty} 
\PP\left(\sup_{\ell = 0,\cdots,\ell_{n}} \left|\frac{Y_{n-\ell}}{a_{n}} - \alpha^{\ell}\right| >  \varepsilon \Big| Y_{n}\geq a_{n}, \ B_{k_{n}}\right).
\end{multline*}
Conditionally on $B_{k_{n}}$,  we can write
$$\frac{Y_{n-\ell}}{a_{n}}=\frac{\alpha^{k_{n}-\ell+1}}{a_{n}}Y_{n-k_{n}-1} + \frac{1}{a_{n}}\sum\limits_{k=0}^{k_{n}-\ell} \alpha^{k}\ep_{n-k-\ell} $$
 where $Y_{n-k_{n}-1}$ and $(\ep_{n-k-\ell} : k=0,\ldots, k_{n}-\ell)$ are still independent  and $(\ep_{n-k-\ell} : k=0,\ldots, k_{n}-\ell)$ are distributed as standard gaussian random variables.  Then,
\begin{equation}\label{eq:termdec}
 \frac{Y_{n-\ell}}{a_{n}} - \alpha^{\ell} = \frac{\alpha^{k_{n}-\ell+1}}{a_{n}}Y_{n-\ell_{n}-1}  +\sum_{k=0}^{k_{n}-\ell} \left(\frac{\alpha^{k}}{a_{n}}\ep_{n-k-\ell} - \alpha^{\ell}(1-\alpha^{2})\alpha^{2k}\right)+\eta_{n,\ell},
\end{equation}
where
\begin{equation}\label{eq:term1}
\eta_{n,\ell}= \alpha^{\ell}(1-\alpha^{2}) \sum_{k=0}^{k_{n}-\ell} \alpha^{2k} - \alpha^{\ell}, \qquad \sup_{\ell \in\{ 0,\ldots,\ell_{n}\}} \vert\eta_{n,\ell}\vert\stackrel{n\rightarrow \infty}{\longrightarrow}0.
\end{equation}
Next, for all $\ell \in\{0,\cdots,\ell_{n}\}$,  \eqref{eq:upyn} yields
\begin{equation*}
\PP\left(\left|\frac{\alpha^{k_{n}-\ell+1}}{a_{n}}Y_{n-k_{n}-1}\right| \geq \ep\right)\leq \PP\left(\left|Y_{n-k_{n}-1}\right| \geq \frac{a_{n}}{\alpha^{k_n-\ell_n}}\ep\right) \lesssim \frac{\alpha^{k_n-\ell_n}}{ a_{n}} \exp\left(-\frac{(1 - \alpha^{2})\ep^{2} \alpha^{-2(k_n-\ell_n+1)} a_{n}^{2}}{2}\right).
\end{equation*}
Dividing now by $\PP\left(Y_{n} \geq a_{n}\right)$ and using again the lowerbound  \eqref{eq:loyn} with $k_n$ yields
\begin{eqnarray*}
&&\sup_{\ell\in\{0,\cdots,\ell_{n}\}} \PP\left(\left|\frac{\alpha^{k_{n}-\ell+1}}{a_{n}} Y_{n - k_{n} - 1}\right| \geq \ep \Big| Y_{n} \geq a_{n}, B_{k_{n}}\right)\\
&& \qquad \qquad \lesssim2^{k_{n}}\alpha^{k_n-\ell_n} \exp\left(-\frac{(1-\alpha^{2})a_{n}^{2}}{2}\left(\ep^{2}\alpha^{-2(k_n-\ell_n+1)}- \frac{1}{1-(\alpha^{2})^{k_{n}+1}}\right)\right) \\
&& \qquad \qquad =\smallO{1/\ell_n}
\end{eqnarray*}
since $k_n-\ell_n\rightarrow \infty$ and $k_n=\smallO{\log{a_{n}}}$.  Using $\mathbb P(\sup_i X_i\geq a)\leq \sum_i \mathbb P( X_i\geq a)$ and the fact that  $Y_{n-k_n-1}$ is independent of $B_{k_{n}}$, it ensures that
\begin{equation}\label{eq:term2}
\lim\limits_{n\rightarrow +\infty}\PP\left(\sup_{\ell\in\{0,\ldots,\ell_{n}\}}\left|\frac{\alpha^{k_{n}-\ell+1}}{a_{n}}Y_{n-k_{n}-1}\right| \geq \ep \Big| Y_{n} \geq a_{n}, B_{k_{n}}\right) = 0.
\end{equation}
Finally, to control the remaining term  in \eqref{eq:termdec}, we prove that
\begin{equation}\label{eq:term3}
\lim_{n\rightarrow +\infty} \PP\Bigg(\sup_{\ell\in\{0,\ldots,\ell_{n}\}} \left|\sum_{k=0}^{k_{n}-\ell} \left(\frac{\alpha^{k}}{a_{n}}\ep_{n-k-\ell} - \alpha^{\ell}(1-\alpha^{2})\alpha^{2k}\right)\right| \geq \ep \Big| Y_{n} \geq a_{n}, B_{k_{n}}\Bigg) = 0.
\end{equation}
It amounts to solve a minimisation problem   for the cost of a  trajectory $(\ep_{n-k-\ell} : \ell\in\{0,\ldots,\ell_{n}\})$, under the constraint $Y_{n} \geq a_{n}$. This problem can be solved explicitly  and we   actually prove in Appendix that for any $k\in \{0,\ldots, k_{n}\}$ \begin{equation}\label{bornegaussienne}
\PP\left(\left|\frac{\alpha^{k}}{a_{n}}\ep_{n-k} - (1 - \alpha^{2})\alpha^{2k} \right| > \ep \Big| Y_{n}\geq a_{n}, B_{k_{n}}\right) \lesssim \exp\left(- C\frac{ a_{n}^{2}}{2}\right),  
\end{equation}
for some positive constant $C$, which yields directly $\eqref{eq:term3}$.

Combining the last estimates ensures that
$$
\lim\limits_{n\rightarrow +\infty}\PP\left(\sup_{\ell\in\{0,\cdots,\ell_{n}\}}\left|\frac{1}{a_{n}}Y_{n-\ell}-\alpha^{\ell}\right| \geq \ep \Big| Y_{n} \geq a_{n}, B_{k_{n}}\right) = 0,
$$
which ends the proof recalling \eqref{partieenv}.
\end{proof}

\subsection{Large deviations  in the weakly stable case}\label{sec:deviation-weakly-stable}

In this regime, $\alpha\beta<1$ and
$$\beta\leq 1<\alpha.$$
 Then there exists a unique $\kappa>0$  such that
$$\alpha^{\kappa}+\beta^{\kappa}=2.$$
We also set $\kappa_0>0$ and $\gamma\in (0,1)$,  respectively  defined by 
$$ \log(\alpha)\alpha^{\kappa_0}+ \log(\beta)\beta^{\kappa_0}=0, \qquad \gamma :=( \alpha^{\kappa_{0}}+\beta^{\kappa_0})/2 =\inf_{s\geq 0}  \{(\alpha^{s}+\beta^{s})/2\}.$$
We introduce now the random walk $(S_k)_{k\geq 1}$ defined by $S_1=0$ and
$$S_{k} = \sum_{\ell=1}^{k-1} \log\theta_{\ell} \quad (k\geq 2).$$
Recalling from the introduction that $A_n^*=A(\theta_n,\ldots,\theta_1)$, we can  write
$$(A_n^*)^2= \sum_{k=1}^{n}
\exp\left(2S_{k}\right).$$ 
We restrict here the study  to the case when $\{A_n^* \geq a_n\}$ is comparable to $\{A_{\infty}^*\geq a_n\}$. In the proof below,  it corresponds to the fact that  the supremum 
$$\overline{S} = \sup_{k\geq 1} S_{k}$$
 is reached before time $n$  when it is larger than $a_n$.   It forces $a_n$  to be not too large in the statement below.
\begin{theorem}\label{thm:large_deviation2}
Let $(a_{n},n\in\NN)$ be a sequence of real numbers  such that 
$a_{n}^{\kappa}.n^{-3/2}.\gamma^n\rightarrow 0$ as $n\rightarrow \infty$. Then, 
\begin{equation}\label{eq:asymp-equiv} 
\PP\left(Y_{n} \geq a_{n}\right) \asymp \PP\left(A_{n}\geq a_{n}\right) \asymp a_{n}^{-\kappa}. \end{equation}
\end{theorem}
 When the supremum is reached latter, additional work is needed and  we expect that other equivalents can be proved using   a change of probability.
\begin{proof}
First we determine the asymptotic behavior of
$\PP\left(A_{n}\geq a_{n}\right)$ using the result of Kesten  \cite{Kesten73} which guarantees that
$(A_{\infty}^*)^2=\sum_{k=1}^{\infty}
\exp\left(2S_{k}\right)$ is comparable to $\exp(2\overline{S})$. More precisely, we first observe 
that
\begin{equation*}
\EE\left[(\theta_{1}^{2})^{\kappa/2}\right] = 1 \quad \text{and} \quad \EE[(\theta_{1}^{2})^{\kappa/2}\log^{+}\theta_{1}^{2}] < \infty.
\end{equation*}
Besides the following equalities hold:
\begin{equation*}
2 S_{k} = \sum_{\ell = 1}^{k-1} \log\theta_{\ell}^{2} \quad \text{and} \quad 2\overline{S} = \sup_{k} \sum_{\ell = 1}^{k-1} \log\theta_{\ell}^{2}.
\end{equation*}
Note that for any $a>0$, $\PP\left(A_{n}^{*} \geq a\right) \leq \PP\left((A_{\infty}^{*})^{2} \geq a^{2}\right).$ Now, using the results of Section XI.6 in \cite{Feller71} and Section 1 in \cite{Kesten73} (see also Section 2 in \cite{ESZ09}), we get, for some positive constant $C$,
\begin{equation*}
\PP\left((A_{\infty}^{*})^{2} \geq a^{2}\right) \sim_{a \rightarrow \infty} C \times \PP\left(\exp\left(2 \overline{S}\right) \geq a^{2} \right) = C \times \PP\left(\exp\left(\overline{S}\right) \geq a\right).
\end{equation*}  
Moreover, we know, see for e.g. Section XI.6 in \cite{Feller71}, that
\begin{equation*}\label{equivmax}
\PP\left(\exp\left(\overline{S}\right) \geq a \right)\sim_{a\rightarrow \infty} c_0a^{-\kappa},
\end{equation*}
where $c_0$ is positive constant. We  conclude  that 
\begin{equation*}
\PP\left(A_{n}^{*} \geq a_{n}\right) \lesssim \PP\left(\exp\left(\overline{S}\right) \geq a_{n}\right) \sim c_0 a_{n}^{-\kappa} \quad \text{when $n \rightarrow \infty$}.
\end{equation*}
 Besides, for any $x\in (0,1]$,
\begin{equation*}
\PP\left(A_{n}^* \geq a_{n}x\right)\geq \PP\left(\sum_{k=1}^{n} \exp\left(2S_{k}\right)  
\geq
a_{n}^{2}x^2, \tau_{\overline{S}} \leq n\right) \geq 
\PP\left(\exp\left(\overline{S}\right) \geq a_{n}x, \tau_{\overline{S}}
\leq n\right),
\end{equation*}
where $\tau_{\overline{S}} := \inf\{k>0; S_{k} = \overline{S}\}$ is the
first hitting time of the maximum of the random walk $(S_{k}, k>0)$.
Adding from Theorem 15,
Chapter 4 of \cite{Borovkov76} that
\begin{equation*}
\PP\left(\tau_{\overline{S}} > n\right) \lesssim n^{-3/2}\gamma^{n},
\end{equation*}
and using  $n^{-3/2}\gamma^{n} = \smallO{a_{n}^{-\kappa}}$, we get that 
 $\PP\left(\tau_{\overline{S}} > n\right)$ is negligible compared to $\PP\left(\exp\left(\overline{S}\right) \geq
a_{n}\right)$ and
$$\PP\left(A_{n}^* \geq a_{n}x\right)\gtrsim
\PP\left(\exp\left(\overline{S}\right) \geq a_{n}x\right)
$$
uniformly for  $n\geq 0$ and $x\in (0,1]$. \\

Combining these estimates and recalling that the symmetry  in law of $(\theta_1,\ldots,\theta_n)$ ensures that $A_n$ is distributed as $A_n^*$,  we get
\begin{equation}\label{major}
\PP\left(A_{n} \geq a_{n}x\right)=\PP\left(A_{n}^* \geq a_{n}x\right) \asymp \PP\left(\exp\left(\overline{S}\right) \geq a_nx\right)
\end{equation}
uniformly for $x\in (0,1]$ and
\begin{equation}\label{major2}
\forall x\in (0,1], \quad \PP\left(\exp\left(\overline{S}\right) \geq a_nx\right) \sim_{n\rightarrow \infty}c_0(xa_{n})^{-\kappa}; \qquad \sup_{x\in (0,1], n\geq 0} \frac{ \PP\left(\exp\left(\overline{S}\right) \geq a_nx\right)}{(xa_{n})^{-\kappa}}<\infty.
\end{equation}
It ensures in particular that $ \PP\left(A_{n}\geq a_{n}\right) \asymp a_{n}^{-\kappa}$ and ends the proof of the second part of the statement. \\

We prove now  that $\PP\left(Y_{n} \geq a_{n}\right) \asymp a_{n}^{-\kappa}.$ 
Setting  $W_{n} = A_{n}/a_{n}$ and recalling that $\mathcal L(Y_n\vert \mathcal F_{n-1})=\mathcal N(0,a_n^2W_n^2)$, we first consider 
\begin{align*}
\PP\left(Y_{n} \geq a_{n}, W_{n} < 1\right) &\sim \EE\left[\indicatrice_{\{W_{n} < 1\}}W_{n}\exp\left(- \frac{1}{2W_{n}^{2}}\right)\right]& \\ &= \int_{0}^{1}\left(1 + \frac{1}{x^{2}}\right) \exp\left(- \frac{1}{2x^{2}}\right) \PP(A_{n} \geq a_nx)dx.& 
\end{align*}
Besides, \eqref{major2} 
and the fact that   $\left(1 + \frac{1}{x^{2}}\right)\exp(- \frac{1}{2x^{2}})x^{-\kappa}$ is integrable on $(0,1]$ ensure  by bounded convergence that  
$$
a_n^{\kappa}\int_{0}^{1}\left(1 + \frac{1}{x^{2}}\right) \exp\left(- \frac{1}{2x^{2}}\right) \PP(\exp(\overline{S}) \geq a_nx)dx\stackrel{n\rightarrow \infty}{\longrightarrow} c \in (0,\infty).
$$  
Using the uniform estimate in \eqref{major} we obtain from the two previous displays that
$$ \PP\left(Y_{n} \geq a_{n}, W_{n} < 1\right)\asymp a_n^{-\kappa} .$$
Next, using again $\mathcal L(Y_n\vert \mathcal F_{n-1})=\mathcal N(0,a_n^2W_n^2)$ and the fact that $\PP(\mathcal N(0,a_n^2x^2)\geq a_n) \in [\PP(\mathcal N(0,1)\geq 1), 1]$ for $x\geq 1$, we get
\begin{align*}
\PP\left(Y_{n} \geq a_{n}, W_{n} \geq 1\right) &\asymp \PP(W_n\geq 1) =\PP(A_n\geq a_n)  \asymp a_{n}^{-\kappa},
\end{align*}
which ends the proof of the theorem.
\end{proof}

\section{Local densities for bifurcating autoregressive processes}

We  can now study the autoregressive process $$Z_n=\sum_{i\in \mathbb G_n} \delta_{X_i}$$ defined in Introduction. We recall that  $\TT$ is the regular binary tree describing the underlying population and $\mathbb{G}_{n}=\{2^{n},2^{n}+1,\cdots,2^{n+1}-1\}$ is the generation $n$.
\subsection{The strictly stable case}$\,$

\begin{theorem}\label{thm:local_density} If $\beta\leq \alpha<1$, then for any $x\in [0,\sqrt{2\log(2)/(1-\alpha^2)})$,
\begin{equation*}
\lim_{n\rightarrow+\infty} \frac{1}{n} \log
Z_{n}([x\sqrt{n};+\infty)) = \log2 - x^2(1-\alpha^{2})/2 \quad \text{a.s}.
\end{equation*}
\end{theorem}
The proof relies on the large deviation results of the previous section and a law of large number principle, which is in the same spirit as  proofs for  branching random walk.  The results of the previous section also tells us that the bulk of   $Z_{n}([x\sqrt{n};+\infty))$ corresponds to individuals coming from the first daughter in the last divisions  and that their value   has undergone  a geometric deviation from their stable distribution in these last divisions.\\

The proof of the upper bound is a classical Borel Cantelli argument using the first moment.
\begin{proof}[Proof of the upper bound of Theorem \ref{thm:local_density}]
Let $\ep>0$ and write
$\alpha(x)= \log2 - x^2(1-\alpha^{2})/2$. By  Markov inequality we have
$$
\PP\left(\frac{1}{n}\log Z_{n}\left([x\sqrt{n};+\infty)\right) \geq \alpha(x)+\ep\right)
  \leq
\exp\left(-n\left(\alpha(x) +
\ep\right)\right) \times
\EE\left[Z_{n}\left([x\sqrt{n};+\infty)\right)\right].$$
Besides, the  many-to-one formula \eqref{MTOf} and \eqref{eq:upyn} yield \begin{equation*}
\EE\left[Z_{n}\left([x\sqrt{n};+\infty)\right)\right] =
2^{n}\PP\left(Y_{n} \geq x\sqrt{n}\right)\lesssim 2^n \exp\left(-\frac{(1-\alpha^{2})x^2n}{2}\right).
\end{equation*}
We obtain
\begin{equation*}
\PP\left(\frac{1}{n}\log Z_{n}\left([\sqrt{n};+\infty)\right) \geq
\alpha(x)+ \ep\right) \leq
\exp\left(-n\ep\right)
\end{equation*}
and  Borel Cantelli Lemma allows us to conclude that
\begin{equation*}
\limsup_{n\rightarrow+\infty} \frac{1}{n} \log
Z_{n}\left([\sqrt{n};+\infty)\right) \leq \alpha(x) \quad \text{a.s.}
\end{equation*}
which ends the proof of the upper bound of Theorem \ref{thm:local_density}.
\end{proof}\

For the lower bound, we need now to prove a law of large number result on the subpopulation producing the local density $[x\sqrt{n},\infty[$ in generation $n$.
Using the results of the previous section, it is achieved by following individuals who undergo 
 a deviation in the last generation in environment $\alpha$.   We first derive from the previous section the following  result on the large deviations of $Y$ and then proceed with the proof of the lower bound.

\begin{lemma} \label{LDderiv} For any $x\geq 0$
\begin{equation}
\liminf_{p\rightarrow \infty} \liminf_{n\rightarrow\infty}  \frac{1}{n+p}\log \PP(Y_p\geq x\sqrt{n+p})=-x^2(1-\alpha^2)/2. 
\end{equation}
\end{lemma}
\begin{proof}
We first observe that by monotonicity of $Y$ with respect to the initial condition and Markov property, $\PP(Y_{n+p}\geq x\sqrt{n+p})\geq \PP(Y_n\geq 0) \PP(Y_p\geq x\sqrt{n+p})= \PP(Y_p\geq x\sqrt{n+p})/2.$
Using \eqref{eq:ldyn}, we get the upper bound 
$$\liminf_{p\rightarrow \infty} \liminf_{n\rightarrow\infty}  \frac{1}{n+p}\log \PP(Y_p\geq x\sqrt{n+p})\leq -x^2(1-\alpha^2)/2.$$
We prove now the converse inequality. We know from Theorem \ref{thm:large_deviation} that conditionally on $\{Y_{n} \geq a_{n}\}$, the process favors the best environment, that is the coefficient $\alpha$, at least in the last time. We then have 
\begin{align*} 
\PP(Y_p\geq x\sqrt{n+p}) &\geq  \PP\left(Y_{p} \geq x\sqrt{n + p}, \theta_{p} = \alpha, \ldots, \theta_{1} = \alpha\right)& \\ &= \frac{1}{2^{p}} \times \PP\left(\sum_{\ell = 0}^{p - 1} \alpha^{\ell} \ep_{p - \ell} \geq x\sqrt{n + p}\right)& \\ &= \frac{1}{2^{p}} \times \PP\left(\Nn\left(0,1\right) \geq \sqrt{\frac{1 - \alpha^{2}}{1 - \alpha^{2p}}}x\sqrt{n + p}\right)& \\ &\gtrsim \frac{1}{2^{p}} \times \frac{1}{x\sqrt{n + p}} \exp\left(- \frac{(1 - \alpha^{2})x^{2}(n + p)}{2(1 - \alpha^{2p})}\right).&
\end{align*}
Now, applying the $\log$ function in both sides of the last inequality, dividing by $n + p$ and letting $n$ and $p$ go to infinity gives the expected lower bound and ends the proof.
 \end{proof}

\begin{proof}[Proof of the the lower bound of Theorem \ref{thm:local_density}.]
We first observe that for any $n,p,a\geq 0$,
$$
Z_{n+p} ([a,\infty))\geq \sum_{u \in \GG_{n} : X_u\geq 0} \#\{ v \in \GG_{n + p} :  \ u \preccurlyeq v, \  X_v\geq a\} \quad \text{a.s}.
$$
Using  the monotonicity of the autoregressive process with respect  to its initial value and  the branching property, 
$$Z_{n+p} ([a,\infty))\geq \sum_{u \in \GG_{n} : X_u\geq 0} X_{p,a}^{(u)} \quad \text{a.s}.$$
where $(X_{p,a}^{(u)})_{u \in \GG_{n}}$  are i.i.d. r.v.   with the same law $\mu_{p,a}$ defined by
\begin{equation*}
\mu_{p,a} = \PP_{\delta_{0}}\left(Z_{p}\left([a;+\infty)\right)\in\cdot\right).
\end{equation*}
Besides, Proposition 28 in \cite{G07} ensures that
$$
\frac{1}{2^n}\#\{u \in \GG_{n} : X_u\geq 0\} \stackrel{n\rightarrow \infty}{\longrightarrow } \int_{\mathbb{R}_{+}} \mu(dx) \quad \text{a.s.},
$$
where $\mu$ is the law of $Y^{*}_{\infty} = \sum_{k = 1}^{\infty} (\prod_{\ell = 1}^{k - 1} \theta_{\ell}) \varepsilon_{k}.$ Actually, note that 
\begin{equation*}
 \int_{\mathbb{R}_{+}} \mu(dx) = \lim_{n \rightarrow \infty} \PP\left(Y_{n} \geq 0\right) = \frac{1}{2}.
\end{equation*}
Setting by now $a=x\sqrt{n+p}$, the many-to-one formula yields
$$u_{n,p}=\overline{\mu_{p,a}}= \EE[X_{p,a}^{(u)}] = \EE[Z_p[x\sqrt{n+p},\infty)] = 2^{p}\PP(Y_p\geq x\sqrt{n+p}). $$
Besides, the large deviation estimate proved above  in Lemma \ref{LDderiv} ensures that
\begin{equation}
\label{adapt}
\liminf_{p\rightarrow \infty} \liminf_{n\rightarrow\infty}  \frac{1}{n+p}\log u_{n,p}=-x^2(1-\alpha^2)/2. 
\end{equation}
Recalling that $x\in [0,\sqrt{2\log(2)/(1-\alpha^2)})$, we consider by now   $p$ large enough such that 
$$\liminf_{n\rightarrow\infty} \frac{1}{n}\log\left(\#\{ u \in \GG_n : X_u \geq 0\}\right)+  \liminf_{n\rightarrow\infty}  \frac{1}{n}\log u_{n,p}= \log(2)+ \liminf_{n\rightarrow\infty}  \frac{1}{n}\log u_{n,p}>0.$$
Writing $(X_{n,i} : i=1 \ldots, N_n)=  (X_{p,a}^{(u)} : u\in \GG_n, X_u\geq 0)$, we are now in position to apply the law of large number given in Proposition \ref{LLNadapted} in Appendix \ref{ALLN} 
and get
$$\liminf_{n\rightarrow \infty} \frac{1}{\#\{ u \in \GG_n : X_u\geq 0\}u_{n,p}}\sum_{u \in \GG_{n} : X_u\geq 0} X_{p,a}^{(u)}=1 \quad \text{a.s}.$$
Then, 
 $$\liminf_{n\rightarrow \infty} \frac{1}{2^{n-1}u_{n,p}} Z_{n+p} ([x\sqrt{n+p},\infty))\geq 1 \quad \text{a.s.}$$
and using \eqref{adapt} we obtain
$$\liminf_{n\rightarrow \infty}  \frac{1}{n} \log
Z_{n}([x\sqrt{n};+\infty)) \geq \log(2)-x^2(1-\alpha^2)/2 \quad \text{a.s.}$$
which ends the proof of the lower bound and  Theorem \ref{thm:local_density}.
\end{proof}

We derive now the asymptotic behavior of the highest value  $M_n=\max\{ X_i : i\in \mathbb G_n\}$.

\begin{corollary}
\label{posextrem}
 $$ \lim_{n\rightarrow+\infty} \frac{1}{\sqrt{n}} M_n=\sqrt{\frac{2 \log(2)}{1-\alpha^{2}}} \qquad \text{a.s}.$$
\end{corollary} 
\begin{proof}
Setting $v_+^2=2 \log(2)/(1-\alpha^{2})+\varepsilon$ for some $\varepsilon>0$, we use the classical estimate
$$\PP(M_n/\sqrt{n}\geq v_+ )\leq \EE(Z_n[\sqrt{n} v_+,\infty) )=2^n\PP(Y_n\geq \sqrt{n}v_+)$$
and recalling \eqref{eq:ldyn},  we get $\sum_n \PP(M_n/\sqrt{n}\geq v_+ )<\infty$. Then  $\limsup_{n\rightarrow \infty} M_n/\sqrt{n} \leq v_+$ a.s. by Borel Cantelli lemma, which proves the upperbound letting $\varepsilon\rightarrow 0$.. \\
Moreover, setting  $v_-=2 \log(2)/(1-\alpha^{2})-\varepsilon$, we observe that  $\{Z([\sqrt{n}v_-,\infty))>0\}\subset \{M_n\geq \sqrt{n}v_-\}$. 
So $$\liminf_{n\rightarrow \infty} M_n/\sqrt{n} \geq v_{-} \quad \text{a.s.}$$
 is a direct consequence of Theorem \ref{thm:local_density}. It  ends the proof.
\end{proof}
\subsection{The weakly stable case}$\,$
We assume now that $\alpha>1$. We recall that $\alpha\beta<1$, so $\beta<1$.  Moreover $\kappa\in (0,\infty)$ is defined by $\alpha^{\kappa}+\beta^{\kappa}=2.$
\begin{theorem}\label{thm:local_density2}
For all $c\in[0,\infty)$ such that $c^{\kappa} < \min(2,1/\gamma)$ and $c>1$, we have
\begin{equation*}
\lim_{n\rightarrow +\infty} \frac{1}{n} \log \#\{i\in \GG_{n} :  X_{i} \geq c^{n}\} = \log(2/c^{\kappa}) \quad \text{in probability}.
\end{equation*}
\end{theorem}
The proof of the upperbound is achieved by a  classical Borel Cantelli argument using the first moment evaluated in the previous section, as for the strictly stable case. The proof of the lowerbound is different and more involved. We need to focus on a subtree producing the large values $[c^n,\infty)$ at time $n$   and characterized by    $A(\theta_1,\ldots,\theta_n) \geq c^n$. We then control correlations through the common ancestor
of nodes and perform  $L^2$ estimates.
\begin{proof}[Proof of the upperbound of Theorem \ref{thm:local_density2}]
Let $\ep>0$. By  Markov inequality we have
\begin{multline*}
\PP\left(\frac{1}{n} \log \# \{i\in \GG_{n} : X_{i} \geq c^{n}\} \geq \log(2/c^{\kappa}) + \ep\right) \leq (c^{\kappa}/2)^n\exp(-n \ep) \times \EE\left[\#\{i\in \GG_{n} :  X_{i} \geq c^{n}\}\right].
\end{multline*}
Besides, using the many-to-one formula and \eqref{eq:asymp-equiv}, 
\begin{equation*}
\EE\left[\#\{i\in \GG_{n} : X_{i} \geq c^{n}\}\right] = 2^{n}\PP\left(Y_{n} \geq c^{n}\right) \lesssim 2^{n} c^{-n\kappa}.
\end{equation*}
From the foregoing we are led to
\begin{equation*}
\PP\left(\frac{1}{n} \log \#\{i\in \GG_{n}  :  X_{i} \geq c^{n}\} \geq \log(2/c^{\kappa}) + \ep\right)  \lesssim \exp\left(-n\ep\right).
\end{equation*}
Finally, the Borel Cantelli Lemma allows us to conclude that
\begin{equation*}
\limsup_{n\rightarrow +\infty} \frac{1}{n} \log \#\{i\in \GG_{n}  :  X_{i} \geq c^{n}\} \leq \log(2/c^{\kappa}) \quad \text{a.s}.
\end{equation*}
which ends the proof of the upper bound.
\end{proof}

Let us turn to the the proof of the lower bound and first give the outline. 
Recall that  $(\theta_{n}, n\geq 1)$ is a sequence of i.i.d. random variables having the same law that $\theta$ and $A_n=A(\theta_1,\ldots,\theta_n)$.
 For any $i\in \GG_n$, we write $(i_{1},\cdots,i_{n})\in\{0,1\}^{n}$  its binary decomposition. This decompositions yields
 the unique path in the  tree from the root 1 to the vertex $i$.
 Setting $\bar{0}=\alpha$ and $\bar{1}=\beta$, we define
 $$A[i] = \sqrt{\sum\limits_{k=1}^{n} \prod\limits_{\ell = k+1}^{n}\overline{i_{l}}^{2}}$$ and we consider the subset of $\GG_n$ defined by 
  $$T_{n} = \{i\in\GG_{n}  :  A[i] \geq c^{n}\}.$$ Note that we have 
\begin{equation}\label{eq:prob-An}
\PP(A_{n}\geq c^{n}) = \frac{\#T_{n}}{2^{n}} \quad \text{and} \quad \#\{i\in\GG_{n} :  X_{i} \geq c^{n}\} \geq \#\{i\in T_{n} :  X_{i} \geq c^{n}\} \quad \text{a.s}.
\end{equation} 
 We set $$F_{n} = \frac{\# \{i\in T_{n}  :  X_{i} \geq c^{n}\}}{\#T_{n}}.$$ The subtree $T_n$ provides the bulk of $Z_n[c^n, \infty)$. Roughly the proof will follow from the fact that 
$\liminf_{n\rightarrow +\infty} F_{n} > 0 \quad \text{in probability}.$
Indeed,  using  Theorem \ref{thm:large_deviation2},  it will ensure that
\begin{equation}
\label{sketch}
\frac{\#\{i\in\GG_{n}  :  X_{i} \geq c^{n}\}}{2^{n}} \geq \frac{\#\{i\in T_{n} : X_{i} \geq c^{n}\}}{2^{n}} = F_{n} \times \PP(A_{n} \geq c^{n})  \gtrsim (c^{n})^{-\kappa}.
\end{equation}
We actually prove   using $L^2$ computations that   $F_n$ is close to the non degenerated sequence 
 \begin{align} f_{n}=\EE[F_{n}] &= \frac{1}{\# T_{n}} \EE\left[\sum\limits_{i\in T_{n}} \mathbf{1}_{\{X_{i} \geq c^{n}\}}\right] = \frac{1}{\# T_{n}} \EE\left[\sum\limits_{i\in \GG_{n}} \mathbf{1}_{\{X_{i} \geq c^{n}\}} \mathbf{1}_{\{i\in T_{n}\}}\right]& \nonumber \\& = \frac{2^{n}}{\# T_{n}} \times \PP(Y_{n} \geq c^{n}, A_n \geq c^{n})
 =\PP(Y_{n}\geq c^{n} | A_{n}\geq c^{n})
 \geq \PP(\Nn(0,1) \geq 1) > 0, \label{positif}\end{align} where $(Y_{n}, n\geq 0)$ is the autoregressive process with random coefficients defined in (\ref{eq:random-ar}) and we use the many to one formula coupling the process and its environment.\\
 For that purpose, for each $i \in \GG_{n}$, we set $$Z_{n}(i) = \mathbf{1}_{\{X_{i} \geq c^{n}\}} - \PP(X_{i} \geq c^{n}).$$ 
 Thus $F_{n} - f_{n} = \frac{1}{\# T_{n}} \sum\limits_{i\in T_{n}} Z_{n}(i)$ and
\begin{eqnarray}
 \EE\left[(F_{n}-f_{n})^{2}\right] &=&
\label{eq:Fn-fn1}
\frac{1}{\# T_{n}^{2}} \sum\limits_{i\in T_{n}} \EE\left[Z_{n}(i)^{2}\right] + \frac{1}{\# T_{n}^{2}}\sum\limits_{\substack{(i,j)\in T_{n}^{2}\\ i\neq j}} \EE\left[Z({i})Z({j})\right].\end{eqnarray} 
For the first term of the right hand side of (\ref{eq:Fn-fn1}), we have $$ \frac{1}{\# T_{n}^{2}} \sum\limits_{i\in T_{n}} \EE\left[Z_{n}(i)^{2}\right] \leq \frac{1}{\# T_{n}^{2}}\sum\limits_{i\in T_{n}} \EE\left[\mathbf{1}_{\{X_{i}\geq c^{n}\}}\right] = \frac{\PP(Y_{n}\geq c^{n}|A_{n}\geq c^{n})}{\# T_{n}} \nlim 0,$$  since $\# T_{n}\rightarrow \infty$ as $n\rightarrow\infty$. Let us deal with the second term of (\ref{eq:Fn-fn1}). For all $p\in\{0,\cdots,n - 1\}$, we denote by $\Tt_{p}^{(n)}$ the set of all the individuals of the generation $p$ who are ancestors of at least one individual in the sub-population $T_{n}$. For each $i\in\GG_{p}$, we denote by $T_{n}(i)$ the set of individuals belonging to $T_{n}$ who are descendants of $i$. 
 Writing  $i \preccurlyeq j$ when $i$ is an ancestor of $j$, it means that 
 $$\Tt_{p}^{(n)} = \{i\in \GG_{p}: \exists  j\in T_{n} \text{ such that } i \preccurlyeq j\} \quad \text{and} \quad T_{n}(i) = \{j\in T_{n}: i\preccurlyeq j\}.$$ 
Writing    $i\wedge j$ the most recent common ancestor of two individuals $i$ and $j$ and gathering the couples in function of their most recent common ancestor, we obtain 
\begin{eqnarray}
&&\frac{1}{\# T_{n}^{2}}\sum\limits_{\substack{(i,j)\in T_{n}^{2}\\ i\neq j}} \EE\left[Z_{n}(i)Z_{n}(j)\right] \nonumber \\
&&\qquad =\frac{1}{\# T_{n}^{2}}\sum\limits_{p=0}^{n-1} \sum\limits_{u\in\Tt_{p}^{(n)}} \sum\limits_{\substack{(i,j)\in T_{n}^{2}\\ i\wedge j = u}} \EE\left[Z_{n}(i)Z_{n}(j)\right] \nonumber \\
\label{decorr} 
&&\qquad = \frac{1}{\# T_{n}^{2}} \sum\limits_{p=0}^{n-1}\sum\limits_{u\in\Tt_{p}^{(n)}} 2\EE\left[\EE\left[\sum\limits_{i\in T_{n}(2u)}Z_{n}(i)\Big|X_{2u}\right] \times \EE\left[\sum\limits_{j\in T_{n}(2u+1)}Z_{n}(j)\Big|X_{2u+1}\right]  \right].\end{eqnarray} 
For all $p$ fixed in $\{0,\cdots,n-1\}$ and  $u\in\GG_p$, we set
 $A_{n}^{p,u}= A(\overline{u_{1}},\ldots,\overline{u_{p}}, \theta_{p+1},\ldots, \theta_n)$ and for any $x\in \mathbb R$,
 $$P_n^p(x,u)=\PP\left(Y_{n} \geq c^{n}\Big| Y_{p} = x, \,A_n^{p,u} \geq c^{n} \right)$$
 and by convention $P_n^p(x,u)=0$ if $\PP (A_n^{p,u} \geq c^{n})=0$. Besides
the many-to-one formula ensures that
\begin{equation}\EE\left[\sum\limits_{i\in T_{n}(u)}Z_{n}(i)\Big|X_{u}\right] =
\#T_{n}(u)\left(P^{p}_{n}(X_{u},u) - \EE\left[P^{p}_{n}(X_{u},u) \right]\right). \label{decorr2}
\end{equation} 
The rest of the proof of the theorem relies on the two following lemmas.
The first one ensures that an  ergodic property holds on lineages, which allows to forget the beginning of the trajectory.
\begin{lemma}\label{lem:p-fix}
For any $p\in \mathbb N$  and  any $u \in \GG_{p}$, we have
\begin{equation*}
P_n^p(X_{u},u) - \EE\left[P_n^p(X_{u},u) \right] \nlim 0 \quad \text{in probability}.
\end{equation*}
\end{lemma}
\begin{proof}[Proof of Lemma \ref{lem:p-fix}]$\,$
Let $p$ be a fixed natural integer and let $u \in \GG_{p}$. First we will prove that conditionally on $\{A_{n} \geq c^{n}\}$, $\prod\limits_{k=1}^{n} \theta_{k}/c^{n}$ converges to 0 in probability. For that purpose we will show that for all $\ep > 0$, we have
\begin{equation}\label{o-theta-k}
\PP\left(\prod\limits_{k=1}^{n} \theta_{k} \geq \ep c^{n}\right) = \smallO{\PP\left(A_{n} \geq c^{n}\right)}.
\end{equation}
We set $$S_{k} = \sum\limits_{\ell=1}^{k} \log(\theta_{\ell}); \quad \overline{S} = \sup\left\{ S_{k}, k \geq 1 \right\} \quad \text{and} \quad \tau_{\overline{S}} = \inf\left\{k\geq 1, S_{k} = \overline{S}\right\}.$$ Let $\xi \in (0,1)$ and $\delta > 0$ such that $\log(\beta)\leq -\delta < \EE[\log(\theta_{1})]<0$. We have
\begin{eqnarray}\label{ineq:theta-k}
\PP\left(\prod\limits_{k=1}^{n} \theta_{k} \geq \ep c^{n}\right) &=&  \PP\left(\prod\limits_{k=1}^{n} \theta_{k} \geq \ep c^{n}, \tau_{\overline{S}} > \xi n\right) + \PP\left(\prod\limits_{k=1}^{n} \theta_{k} \geq \ep c^{n}, \tau_{\overline{S}} \leq \xi n\right)\nonumber  \\ &\leq& \PP\left(\tau_{\overline{S}} > \xi n\right) +  \PP\left(\overline{S}> n(\log(c) + \delta)\right) \nonumber \\ 
&& \quad + \sum\limits_{j=1}^{\xi n} \PP\left(\prod\limits_{k=1}^{n} \theta_{k} \geq \ep c^{n}, \tau_{\overline{S}} = j, S_{j} \leq n(\log(c) + \delta)\right).
\end{eqnarray}
Using the asymptotic of $\tau_{\overline{S}}$ given in Theorem 15,
Chapter 4 of \cite{Borovkov76}, we obtain
\begin{equation}\label{ineq:taubar1}
 \PP\left(\tau_{\overline{S}} > \xi n\right) \lesssim (\xi n)^{-3/2} \gamma^{\xi n}. 
\end{equation}
Next,  from Kesten's results, see  e.g. \cite{Kesten73}, Section 1 or \cite{ESZ09} Section 2,  we have 
\begin{equation}\label{ineq:taubar2}
\PP\left(\overline{S} \geq n(\log(c) + \delta)\right) \sim (c^{n})^{-\kappa}\exp(-\kappa\delta n).
\end{equation} 
Finally, for the last term of the inequality (\ref{ineq:theta-k}) we have
\begin{multline*}
\sum\limits_{j=1}^{\xi n} \PP\left(\prod\limits_{k=1}^{n} \theta_{k} \geq \ep c^{n}, \tau_{\overline{S}} = j, S_{j} \leq n(\log(c) + \delta)\right) \\ \leq \sum\limits_{j=1}^{\xi n} \PP\left( \exp(S_n-S_j) \geq \ep \exp\left(-n\delta\right), \tau_{\overline{S}} = j, S_{j}\geq S_n \geq n\log(c)+\log(\ep)\right) \\ = \sum\limits_{j=1}^{\xi n} \PP\Bigg(\tau_{\overline{S}} = j, S_{j}>n\log(c)+\log(\ep)\Bigg)\PP\left(\exp(S_n-S_j) \geq \ep \exp(-n\delta) \Big| \tau_{\overline{S}} = j\right).
\end{multline*}
We first observe that conditioning on $\tau_{\overline{S}} = j$ make stochastically decrease the random after  time $j$ and
\begin{eqnarray*}
\PP\left(\exp(S_n-S_j) \geq \ep \exp(-n\delta) \Big| \tau_{\overline{S}} = j\right)& \leq &\PP\left(S_n-S_j\geq  -n\delta +\log(\ep) \right) \\
&\lesssim & \PP\left(S_{n-j}\ \geq -n\delta \right) \\
&\leq &
\exp\left(-(n-j)\psi\left(-\frac{n}{n-j}\delta\right)\right),
\end{eqnarray*}
where 
$\psi$ is the rate function associated to $(S_{k}, k\geq 1)$,  see e.g. \cite{DZ98} Chapter 2 for more details.
Since $\psi$ is nonincreasing in $(-\infty,\EE[\log(\theta_{1})]]$, we have for $j\in\{1,\cdots,\xi n\}$,
\begin{equation*}
\exp\left(-(n-j)\psi\left(-\frac{n}{n-j}\delta\right)\right) \lesssim \exp\left(-(1-\xi)n\psi\left(-\delta\right)\right).
\end{equation*}
From the foregoing we obtain
\begin{multline}\label{ineq:taubar3}
\sum\limits_{j=1}^{\xi n} \PP\left(\prod\limits_{k=1}^{n} \theta_{k} \geq \ep c^{n}, \tau_{\overline{S}} = j, S_{j} \leq n(\log(c) + \delta)\right) \\ \lesssim \PP\left(S_{\tau_{\overline{S}}} > n\log(c)\right) \times \exp\left(-(1-\xi)n\psi\left(-\delta\right)\right) \lesssim (c^{n})^{-\kappa} \times \exp\left(-(1-\xi)n\psi\left(-\delta\right)\right).
\end{multline}
From (\ref{ineq:taubar1})-(\ref{ineq:taubar3}) we conclude that (\ref{o-theta-k}) holds and then that
\begin{equation*}
\PP\left(\prod\limits_{k=1}^{n} \theta_{k} \geq \ep c^{n} \Big| A_{n} \geq c^{n}\right) \nlim 0.
\end{equation*}
Similarly, it is easy to see that for $p\in \mathbb N$ and $u\in\GG_{p}$ and $x\in \RR$, we also have
\begin{equation}\label{eq:thetakXu}
\PP\left(\Bigg|x \, \frac{\prod_{k=p+1}^{n}\theta_{k}}{c^{n}} \Bigg|\geq \ep \Big| A_{n}^{p,u} \geq c^{n} \right) \nlim 0.
\end{equation}
We observe now that
\begin{equation}\label{eq:Yn-p}
P_n^p(x,u) = \PP\left(x\, \frac{\prod_{k=p+1}^{n} \theta_{k}}{c^n}  + \frac{\widetilde{Y}_n}{c^n}\geq 1  \Big| A_n^{p,u} \geq c^{n}  \right),
\end{equation}
where $ \widetilde{Y}_n =\sum_{k=p+1}^{n} \left(\prod_{\ell = k+1}^{n} \theta_{l}\right)\ep_{k}$. 
Combining \eqref{eq:Yn-p}, \eqref{eq:thetakXu} and the fact that $\widetilde{Y}_n/c^n$ is a mixed of Gaussian random variable and the fact that $X_u$ is stochastically bounded (since gaussian),  we get
\begin{equation*}
P_n^p(X_{u},u)  - \PP\left(\widetilde{Y}_{n} \geq c^{n}\Big| A_{n}^{p,u} \geq c^{n}\right)\rightarrow 0 \quad \text{in probability}.
\end{equation*}
By dominated convergence theorem, we obtain that
$\EE\left[P_n^p(X_{u},u) \right]- \PP\left(\widetilde{Y}_{n} \geq c^{n}\Big| A_{n}^{p,u} \geq c^{n}\right)$
goes to $0$ as $n$ tends to infinity. The proof follows by combining these two  limits.
\end{proof}
The next lemma proves that the common ancestor of lineages of the bulk $Z_n[c^n,\infty)$ are at the beginning of the tree, thus providing the decorrelation needed for a law of large number.
\begin{lemma}\label{lem:p-large}
\begin{equation*}
\sup_{n > p}\left\{ \frac{\sum\limits_{k=p}^{n-1}\sum\limits_{u\in\Tt_{k}^{(n)}} \#T_{n}(2u) \,  \#T_{n}(2u+1)}{\# T_{n}^{2}}\right\} \plim 0.
\end{equation*}
\end{lemma} 
\begin{proof}[Proof of  Lemma \ref{lem:p-large}]$\,$
Let $I_{n}$ and $J_{n}$ be two independent indices uniformly drawn from $T_{n}$ and independent of $(X_{i}, i\in \TT)$. Let $p < n$. We have
\begin{equation*}
\PP(p\leq |I_{n}\wedge J_{n}| <n) = \frac{1}{\# T_{n}^{2}} \sum\limits_{(i,j)\in T_{n}^{2}, i\ne j } \indicatrice_{\{|i\wedge j|\geq p\}}. 
\end{equation*}
Next, gathering pairs in function of their most recent common ancestor we obtain
\begin{equation*}
\PP(p\leq |I_{n}\wedge J_{n}| <n ) = \frac{1}{\# T_{n}^{2}} \sum\limits_{(i,j)\in T_{n}^{2}, i\ne j} \indicatrice_{\{|i\wedge j|\geq p\}} = \frac{1}{\# T_{n}^{2}}\sum\limits_{k=p}^{n-1}\sum\limits_{u\in\Tt_{k}^{(n)}} \# T_{n}(2u). \#T_{n}(2u+1).
\end{equation*}
Since for any $u\in \Tt_{p}^{(n)}$, $$\PP\left(p\leq |I_{n}\wedge J_{n}|<n \left|\right. u\preccurlyeq I_{n}\right) = \frac{\# T_{n}(u)}{\# T_{n}},$$ then
\begin{equation*}\label{eq:boundInJn}
\PP\left(p\leq |I_{n}\wedge J_{n}|<n\right) \leq \frac{\max\left\{\# T_{n}(u): u\in\Tt_{p}^{(n)}\right\}}{\# T_{n}}.
\end{equation*}
Moreover, since $\alpha > \beta$, we easily observe that
\begin{equation*}
\max\left\{\# T_{n}(u): u\in\Tt_{p}^{(n)}\right\} = \# T_{n}(2^{p}).
\end{equation*}
Thus we need to check that
\begin{equation}\label{eq:limT2p}
\frac{\sup\limits_{n > p} \left\{\# T_{n}(2^{p}) \right\}}{\# T_{n}} \plim 0. 
\end{equation}
For that purpose, we set for $p < n$,
\begin{equation*}
N_1(i)= \# \left\{k: i_{k} = 1 \right\}, \qquad B_{n}^{p} = \left\{i\in T_{n}(2^{p}): N_1(i) < \ep(n-p)\right\},
\end{equation*}
where $\ep > 0$ is fixed. We recall that for any $i\in\GG_{n}$, $(i_{1},\cdots,i_{n})\in \{0,1\}^n$ is the binary decomposition providing the unique path in the binary tree from the root to $i$. In word, $B_{n}^{p}$ is the set of individuals of $T_{n}(2^{p})$ whose binary decomposition has less than $\ep(n-p)$ ``1''. We first check that there exists $\ep > 0$ such that
\begin{equation}\label{eq:limBnp}
\frac{\#B_{n}^{p}}{\# T_{n}} \plim 0.
\end{equation}
It means  that the number of extremal individuals created by the prolific individual ``$2^{p}$'' is negligible with respect to the total number of extremal individuals when $p$ becomes large. Now let $\psi$ be rate function associated to the simple random walk (with step 0 and 1 with probability 1/2). For $n>p$ large enough, we have
\begin{equation*}
\frac{\# B_{n}^{p}}{2^{n-p}} = \frac{\# \left\{i\in \GG_{n-p}: N_1(i)<\ep(n-p)\right\}}{2^{n-p}} \leq \exp\left(-(n-p)\psi(\ep)\right). 
\end{equation*}
 Recalling that $\#T_{n} \sim (2/c^{\kappa})^{n}$, there exists $c_0>0$ such that
$$
\frac{\# B_n^p }{ \# T_n} \leq c_0\left(\frac{\exp\left(\psi(\varepsilon)\right)}{2}\right)^p \left(\frac{c^{\kappa}}{\exp\left(\psi(\varepsilon)\right)}\right)^n.
$$
 Now we recall that  $c^{\kappa}<2$ and $\psi(0)=\log(2)$. So we can choose $\varepsilon >0$ such that $c^{\kappa}<\exp\left(\psi(\varepsilon)\right)$. Then
 $$\frac{\# B_n^p}{ \#T_n} \leq  c_0\left(\frac{\exp\left(\psi(\varepsilon)\right)}{2}\right)^p .\left(\frac{c^{\kappa}}{\exp\left(\psi(\varepsilon)\right)}\right)^p\leq \left(\frac{c^{\kappa}}{2}\right)^p$$
for $n > p $. The last term tends to $0$ as $p$ tends to infinity, which yields \eqref{eq:limBnp}. 

Now using \eqref{eq:limBnp}, a sufficient condition for \eqref{eq:limT2p} is
\begin{equation}\label{condsuffis}
\frac{\#  \, T_{n}(2^{p})\setminus B_{n}^{p}}{\#T_{n}} \plim 0,
\end{equation}
where $T_{n}(2^{p})\setminus B_{n}^{p}$ is the set of individuals of $T_{n}(2^{p})$ whose the number of components which are equal to 1 is greater than $\ep(n-p)$. We prove \eqref{condsuffis} as follows. For $i\in \GG_n$ and $1\leq a<b\leq n$, we define $\tau_{ab} i\in \GG_n$ as the label obtained from $i$ by permutation of $i_a$ and $i_b$ in the binary decomposition (or path in the binary tree). We also introduce
\begin{equation*}
\mathfrak I_n^k =\left\{ u \in T_n(2^{p}) : N_1(u)=k\right\}
\end{equation*}
the set of individuals in $T_n$, whose ancestor in generation $p$ is $2^{p}$, and which contains exactly $k$ components equal to $1$.\\ Let us first fix $a$ with $1\leq a\leq p$ and consider the set
\begin{equation*}
\mathfrak I_n^k(a) =\{ \tau_{a b}u :   u \in \mathfrak I_n^k, \ p+1\leq b\leq n, \  u_b=1\}.
\end{equation*}
First $\mathfrak I_n^k(a)\cap T_n(2^{p}) = \emptyset $ and $\mathfrak I_n^k(a)\subset T_n$ since  
by definition of $A_n$ and using   $\beta<\alpha$, $A_n(\tau_{ab}u )\geq A_n(u)$ for any $u$ such that $u_a=0$ and $u_b=1$. 
Moreover for any $a\ne a'$ and $k\ne k'$, $\mathfrak I_n^k(a)\cap \mathfrak I_n^{k'}(a')=\emptyset$, so
\begin{equation*}
\# \, T_n\setminus T_n(2^{p}) \,  \geq  \, \# \bigcup_{\substack {k \geq (n-p)\varepsilon \\  1\leq a \leq p}} \mathfrak I_n^k(a)=\sum_{\substack{k \geq (n-p)\varepsilon, \\  1\leq a \leq p}} \#\mathfrak{I}_n^k(a).
\end{equation*}
Besides,
\begin{equation*}
\# \mathfrak I_n^k(a) = \frac{k}{n-p-k} \# \mathfrak I_n^k,
\end{equation*}
since 
\begin{itemize}
\item for any $u \in \mathfrak I_n^k$, we have $k$ choices for $b$ such that $u_b=1$;
\item for any $v\in  \mathfrak I_n^k(a)$, we have $n-p-k$ choices for $0$ in the $n-p$ last components of $v$ that you can change for $1$ : 
$\# \left\{  (u,b) :  u  \in \mathfrak I_n^k, \ p+1\leq b\leq n, \  u_b=1,  \tau_{a b}u =v \right\}=(n-p-k).$
\end{itemize}
Putting the three last facts together, we get
\begin{equation*}
\# \, T_n\setminus T_n(2^{p}) \geq \sum_{ \substack{(n-p)\varepsilon\leq k< n-p, \\  1\leq a \leq p}}\frac{k}{n-p-k} \# \mathfrak I_n^k \geq  \varepsilon  \sum_{ \substack{(n-p)\varepsilon\leq k< n-p, \\  1\leq a \leq p}} \# \mathfrak I_n^k.
\end{equation*}
Excluding the path $(0,\cdots,0,1,\cdots,1)$ corresponding to $k=n-p$,   $T_n(2^{p})\setminus B_n^p$  is the union of the disjoint sets $\mathfrak I_n^k$ for $(n-p)\varepsilon\leq k< n-p$. Thus, we get
\begin{equation*}
\#  \, T_n\setminus T_n(2^p) \,  \geq \,
\varepsilon p \times \# \, T_n(2^{p})\setminus B_n^p
\end{equation*}
and
\begin{equation*}
\lim_{p\rightarrow \infty} \sup_{n\geq p} \frac{\#  \, T_n(2^{p})\setminus B_n^p}{\# \, T_n\setminus T_n(2^{p})}=0.
\end{equation*}
This proves \eqref{condsuffis} and ends the proof of  Lemma \ref{lem:p-large}.
\end{proof}
We can now proceed with the proof of Theorem \ref{thm:local_density2}.
\begin{proof}[Proof of the lower bound of Theorem \ref{thm:local_density2}]
Now let $\ep > 0$. From Lemma \ref{lem:p-large} let $p_{\ep}$ such that
\begin{equation*}
\frac{1}{\# T_{n}^{2}}\sum\limits_{k=p_{\ep}}^{n-1}\sum\limits_{u\in\Tt_{k}^{(n)}} \#T_{n}(2u)\,  \#T_{n}(2u+1) \leq  \ep.
\end{equation*}
Adding that $\vert Z_n(i)\vert \leq 1$ and splitting the sum for $p$ larger than $p_{\ep}$, we get
$$
\frac{1}{\# T_{n}^{2}}\sum\limits_{\substack{(i,j)\in T_{n}^{2}\\ i\neq j}} \EE\left[Z_n({i})Z_n({j})\right] \leq
R_n+ \ep, 
$$
where
$$R_n=\frac{1}{\# T_{n}^{2}}\sum\limits_{p=0}^{p_{\ep}-1} \sum\limits_{u\in\Tt_{p}^{(n)}} 
\sum\limits_{\substack{(i,j)\in T_{n}^{2}\\ i\wedge j = u}} \EE\left[Z_{n}(i)Z_{n}(j)\right].$$
Besides, recalling \eqref{decorr}-\eqref{decorr2},
\begin{multline*}
R_{n} = \frac{1}{\# T_{n}^{2}}\sum\limits_{p=0}^{p_{\ep}-1}  \sum\limits_{u\in\Tt_{p}^{(n)}} 2 . \# T_{n}(2u). \#T_{n}(2u+1) \\ 
\times [P_n^{p+1}(X_{2u},2u)-\EE(P_n^{p+1}(X_{2u},2u))][P_n^{p+1}(X_{2u+1},2u+1)-\EE(P_n^{p+1}(X_{2u+1},2u+1))].
\end{multline*}
Next, $P_n^p$ is bounded and by dominated convergence theorem,  Lemma \ref{lem:p-fix}  yields $R_{n} \nlim 0.$ We conclude from the foregoing that
\begin{equation*}
\frac{1}{\# T_{n}^{2}} \sum_{\substack{(i,j)\in T_{n}^{2}\\ i\neq j}} \EE\left[Z_{n}(i)Z_{n}(j)\right] \nlim 0
\end{equation*}
and then that
\begin{equation*}
\EE\left[(F_{n} - f_{n})^{2}\right] \nlim 0.
\end{equation*}
Combining this limit with \eqref{sketch} and \eqref{positif}  ends the proof of the lowerbound and Theorem \ref{thm:local_density2}.
\end{proof}

Finally, we derive the asymptotic value of the maximal value $M_n$ in the weakly stable case.
\begin{corollary}
 $$ \lim_{n\rightarrow+\infty} \frac{1}{n} \log(M_n)=\frac{\log(2)}{\kappa} \qquad \text{in probability}.$$
 \end{corollary}
The proof follows  the proof of Corollary \ref{posextrem} in the strictly stable case and is left to the reader. One need now to use the estimate of Theorem \ref{thm:large_deviation2} for the  upper bound and check inequalities in probability for the lower bound.

\section{Appendix}

\subsection{Proof of  \eqref{bornegaussienne} : minimization problem for  large values
in the strictly stable case}$\,$\newline
\label{AppGau}
We first focus on the minimization problem and then use it to get the estimates of probabilities given in \eqref{bornegaussienne}.
We set $\vartheta =1/\alpha^2>1$ for convenience and consider the
space 
$$
H=\{ y\in \mathbb R^{\mathbb N} :  \sum_{k\geq 0} y_k^2\vartheta^k  < \infty\}.
$$
It is easy to see that $H$ is an Hilbert space endowed with the inner product $<x,y>= \sum_{k\geq 0} x_k y_k\vartheta^k$. We write  $ \parallel y \parallel^{2} = \sum_{k\geq 0} y_k^2\vartheta^k$ the associated norm. Thanks to  Cauchy-Schwarz inequality, one can observe that $y\rightarrow \sum_{k\geq 0} y_k$ is well defined and continuous on $H$. For $D\subset H$, we define
\begin{equation}\label{eq:ID}
I(D)= \inf \left\{ \parallel y \parallel^2 : \sum_{j\geq 0} y_j \geq 1, y\in D\right\}.
\end{equation}
and we have the following minimization result for a quadratic form with a linear constraint.
\begin{lemma}\label{lem:opti}  
Writing $v=((1-\alpha^2)\alpha ^{2j})_{j\geq 0}$, 
for all $k\leq n$ and for all $\ep > 0$, we have
\begin{eqnarray*}
&i)&I(H)=\parallel v\parallel^2 =1-\alpha^2;\\
&ii)&   I(\{ y \in H :  \vert y_{k} - v_k\vert \geq \varepsilon\})> I(H);\\
&iii)& I(\{y\in H :  \ \vert y_{k} -v_k\vert \geq \varepsilon, y_i=0 \text{ for } i>n\})\geq  I(\{ y \in H :  \vert y_{k} - v_k\vert \geq \varepsilon\}).
\end{eqnarray*}
\end{lemma}
\begin{proof} 
First, since $\parallel y \parallel^2$ is invariant by change of sign of a coordinate, we get
$$
I(H)= \inf \left\{ \parallel y \parallel^2 : \sum_{j\geq 0} y_j = 1, y\in H\right\}.
$$
We observe that for all $y \in H$ such that $\sum_{j\geq 0} y_j = 1,$ the vector $z=y-v$ belongs to $H$ and satisfies $\sum_{j\geq 0} z_j = 0.$ Then $y=v+z$ and  $<v,z>=0$ and
$$ 
\parallel y \parallel^2=\parallel v \parallel^2 +2<v,z>+\parallel z \parallel^2=\parallel v\parallel^2+\parallel z \parallel^2
$$
which proves $i)$ and $ii)$. Next, $iii)$ is a direct consequence of 
$$
\{y\in H :  \ \vert y-v_k\vert \geq \varepsilon, \, \, \, \, \text{and} \, \, \, \, y_i=0 \, \, \, \, \text{for} \, \, \, \, i>n\}\subset \{ y \in H :  \vert y-v_k\vert \geq \varepsilon\}.
$$
which ends the proof.
\end{proof}
 
We can now deal with the proof of \eqref{bornegaussienne}. We recall that $(k_{n},n\in\NN)$ is a sequence of integers smaller than $n$ and tending to infinity such that $k_n=\smallO {\log(a_n)}$ and 
$
B_{k_{n}}=\{ \theta_{n} = \alpha,\cdots,\theta_{n-k_{n}} = \alpha\}.
$  
For each $n\geq 1$, we  subdivise $[-1,1]$ as a collection of $N_n$ successive disjoint intervals $(I_i^n : 1\leq i\leq N_{n})$ of length (at most) $1/n^{3/2}:$
\begin{equation*}
[-1,1]=\cup_{i=1}^{N_n} I_i^n, \quad  \vert I_{i}^n\vert \leq 1/n^{3/2}, \quad I_{i}^n\cap I_{j}^n=\varnothing \text{ for } i\ne j, \quad  N_n\leq 2(n^{3/2}+1)
\end{equation*}
and either $I_i^n\subset \RR_{+}$ or $I_i^n\subset \RR_{-}$. We set $I_{1}^{n} = [0,1/n^{3/2}]$, $I_{2}^{n} = [-1/n^{3/2},0)$ and $I_{0}^{n} = (-\infty,-1)\cup(1,+\infty)$. Then we can write $\RR = \cup_{i=0}^{N_n} I_i^n.$ For $k=0,\ldots, k_{n}$ and $\ep > 0$, we write $\mathcal I^n(k)$ the set of paths ${\bf i} \in \{0,\ldots,N_{n}\}^{n}$  such that the $k^{th}$ interval, $I_{{\bf i}_{k}}$, does not intersect $[v_{k} - \ep ,
v_{k} + \ep]$ :
$$\mathcal I^n(k) =\{{\bf i} \in \{0, \ldots ,N_n\}^{n} :  I_{{\bf i}_{k}}^n\cap [v_{k} - \ep ,
v_{k} + \ep]=\varnothing \},$$
where we recall that $v_k= (1-\alpha^{2})\alpha^{2k}$.
We define
$$
\Xi_{j,n}=\frac{\alpha^{k_{n} + 1}}{a_{n}} \theta_{n - k_{n} - 1} \ldots \theta_{n - j + 1} \, \ep_{n - j},
$$ 
for  $j\in \{k_{n} + 2, \ldots, n\}$ and  
$
\Xi_{j,n}=\frac{\alpha^{k_{n} + 1}}{a_{n}}  \ep_{n - j},
$
for  $j \in \{0,\ldots,k_{n} + 1\}$. We set 
\begin{equation*}
E_{n}({\bf i}) = \bigcap_{ j = 0}^{n} \left\{\Xi_{j,n} \in I_{{\bf i}_j}^n\right\}.
\end{equation*}
 For all $\ep > 0$, we have  
\begin{equation}\label{eq:cvTkn}
\limsup_{n\rightarrow \infty} \PP\left(\left|\frac{\alpha^{k}}{a_{n}}\ep_{n-k} - v_{k}
\right| > \ep \Big| Y_{n}\geq a_{n}, B_{k_{n}}\right)
\leq \limsup_{n\rightarrow \infty}  \sum_{{\bf i} \in \mathcal{I}^{n}(k)} \PP\left(E_{n}({\bf i}) | Y_{n} \geq a_{n}, B_{k_{n}}\right).
\end{equation}

Let $r \in (0,1)$ fixed. We set $\mathcal{I}_{1}^{n}(k) = \{{\bf i} \in \mathcal{I}^ {n}(k): ({\bf i}_{n^{r}}, \cdots, {\bf i}_{n - 1}) = \{1,2\}^{n-n^{r}}\}$ and
\begin{equation*}\label{eq:Dnk}
\delta_{n} = 1 - \frac{(n - n^{r})}{n^{3/2}} \quad \text{and} \quad D_{n,k} = \left\{y  \in \mathbb{R}^{n}; \quad \sum_{\ell = 0}^{n^{r}-1} y_{\ell}\geq 1 - \delta_{n}; \quad \left|y_{k}-v_{k}\right|\geq \ep\right\}.
\end{equation*}
Then, we have the following bounds.
\begin{lemma}\label{lem:Est-Xi} 
For all ${\bf i} \in \mathcal{I}^{n}_{1}(k)^{c}$,
\begin{equation}\label{eq:claim1}
\PP\left(E_{n}({\bf i}) | Y_{n} \geq a_{n}, B_{n}\right) \lesssim \alpha^{n^{r}} n^{-3/2} u_{n}^{-1} \exp\left(- \frac{a_{n}^{2}}{2}\left(\left(\alpha^{-n^{r}} n^{-3/2}\right)^{2} - u_{k_{n}}^{-2}\right)\right).
\end{equation}
 For  ${\bf i} \in \mathcal{I}^{n}_{1}(k)$,
\begin{equation}\label{eq:claim2}
\PP\left(E_{n}({\bf i}) | Y_{n} \geq a_{n}, B_{k_{n}}\right)  \lesssim (1/2)^{n - n^{r}} (n^{3/2})^{n^{r}} R_n,
\end{equation}
where
\begin{equation*}
R_n = \exp\left( - \frac{a_{n}^{2}}{2}\left(\inf\limits_{y \in D_{n,k}} \sum_{j = 0}^{n^{r} - 1} \alpha^{-2j} y_{j}^{2} - \frac{1 - \alpha^{2}}{1 - \alpha^{2k_{n} + 2}}\right)\right).
\end{equation*}
\end{lemma}
\begin{proof}[Proof of \eqref{eq:claim1}.]
First, note that for all $j \geq n^{r},$ we have
\begin{eqnarray*}
\PP\left(\Xi_{j,n} \geq \frac{1}{n^{3/2}} \, \big\vert \, Y_{n} \geq a_{n}, B_{k_{n}}\right)& \leq &\frac{\PP\left(\Xi_{n^{r},n} \geq \frac{1}{n^{3/2}}\right) \PP\left(B_{k_{n}}\right)}{\PP\left(Y_{n} \geq a_{n}, B_{k_{n}}\right)} \\
&\leq & \frac{\PP\left(\ep_{n - n^{r}} \geq a_{n}\alpha^{-n^{r}}n^{-3/2}\right) \PP\left(B_{k_{n}}\right)}{\PP\left(Y_{n} \geq a_{n}, B_{k_{n}}\right)} \nonumber \\ &  \lesssim & \alpha^{n^{r}} n^{-3/2} u_{k_{n}}^{-1} \exp\left(- \frac{a_{n}^{2}}{2}\left(\left(\alpha^{-n^{r}} n^{-3/2}\right)^{2} - u_{k_{n}}^{-2}\right)\right). \label{eq:BXinrn}
\end{eqnarray*}
In the same way, we have
\begin{equation*}
\PP\left(\Xi_{j,n} \leq -\frac{1}{n^{3/2}} \, \big\vert \, Y_{n} \geq a_{n}, B_{k_{n}}\right) \lesssim \alpha^{n^{r}} n^{-3/2} u_{k_{n}}^{-1} \exp\left(- \frac{a_{n}^{2}}{2}\left(\left(\alpha^{-n^{r}} n^{-3/2}\right)^{2} - u_{k_{n}}^{-2}\right)\right),
\end{equation*}
which implies that
\begin{multline}\label{eq:claim0}
\forall j \geq n^{r}, \quad \PP\left(|\Xi_{j,n}| \geq \frac{1}{n^{3/2}} | Y_{n} \geq a_{n}, B_{k_{n}}\right) \\ \lesssim \alpha^{n^{r}} n^{-3/2} u_{k_{n}}^{-1} \exp\left(- \frac{a_{n}^{2}}{2}\left(\left(\alpha^{-n^{r}} n^{-3/2}\right)^{2} - u_{k_{n}}^{-2}\right)\right).
\end{multline}
Now, for ${\bf i} \in \mathcal{I}^{n}_{1}(k)^{c}$, let $j_{n} \geq n^{r}$ such that ${\bf i}_{j_{n}} \neq -1.$  Using \eqref{eq:claim0}, we get
\begin{align*}
\PP\left(E_{n}({\bf i}) | Y_{n} \geq a_{n}, B_{k_{n}}\right) &\leq \PP\left(\bigcap_{ j = 0}^{n^{r} - 1} \left\{\Xi_{j,n} \in I_{{\bf i}_j}^n\right\} \bigcap \left\{ \Xi_{j_{n},n} \in I^{n}_{{\bf i}_{j_{n}}} \right\} \Big| Y_{n} \geq a_{n}, B_{k_{n}}\right)& \\  &\leq  \PP\left(|\Xi_{j_{n},n}| \geq \frac{1}{n^{3/2}} | Y_{n} \geq a_{n}, B_{k_{n}}\right)& \\ &\lesssim \alpha^{n^{r}} n^{-3/2} u_{k_{n}}^{-1} \exp\left(- \frac{a_{n}^{2}}{2}\left(\left(\alpha^{-n^{r}} n^{-3/2}\right)^{2} - u_{k_{n}}^{-2}\right)\right)&  
\end{align*}
and this ends the proof of \eqref{eq:claim1}.
\end{proof}
\begin{proof}[Proof of \eqref{eq:claim2}]
We set $y_i$ the lower bound of $I_i$ if $I_i\subset \RR_+$ and the upper bound of $I_i$ if $I_i\subset \RR_{-}$. 
We note that for all $j \in \{n^{r}, \cdots, n - 1\},$ we have $\{\Xi_{j,n} \in I_{1}^{n}\} \subset \{\ep_{n-j} \geq 0\}$ and $\{\Xi_{j,n} \in I_{2}^{n}\} \subset \{\ep_{n-j} < 0\}.$ We also note that, for all $j \in \{n^{r},\ldots,n-1\},$ we have $-1/n^{3/2} \leq \Xi_{j,n} \leq 1/n^{3/2}$ and then that $$1 - \frac{(n - n^{r})}{n^{3/2}} \leq 1 - \sum_{j = n^{r}}^{n - 1} \Xi_{j,n} \leq 1 + \frac{(n - n^{r})}{n^{3/2}}.$$ Now since $Y_{n}/a_{n} = \sum_{j = 0}^{n^{r} - 1} \Xi_{j,n} +  \sum_{j = n^{r} }^{n - 1} \Xi_{j,n},$ we get 
\begin{equation}\label{eq:Ipir0}
\PP\left( E_{n}({\bf i}), Y_{n}\geq a_{n}, B_{k_{n}} \right) \leq (1/2)^{n-n^{r}} \PP\left(D_n({\bf i})\right),
\end{equation}
where 
$$D_n({\bf i})=\bigcap_{ j = 0}^{n^{r} - 1} \left\{\Xi_{j,n} \in I_{{\bf i}_j}^n\right\} \bigcap \left\{\sum_{j = 0}^{n^{r} - 1} \Xi_{j,n} \geq 1 - \delta_{n} \right\}\bigcap B_{k_{n}}.$$
Using the convention $0/0 = 0$ and that 
$$
D_n({\bf i})\,  \subset \, \bigcap_{ j = 0}^{n^{r} - 1} \left\{\Xi_{j,n} \in I_{{\bf i}_j}^n\right\} \bigcap  \left\{ \sum_{j = 0}^{n^{r} - 1} (y_{{\bf i}_{j}} + n^{-3/2}) \geq 1 - \delta_{n} \right\} \bigcap B_{k_{n}},$$
we get
\begin{align*}
 \PP\left(D_n({\bf i})\right)
 & \leq \PP\left(B_{k_{n}}\right)\left(\prod_{j = 0}^{n^{r} - 1}\left(\frac{\alpha^{j}}{a_{n}|y_{{\bf i}_{j}}|}\mathbf{1}_{\{y_{{\bf i}_{j}} \neq 0\}} + \frac{1}{2} \mathbf{1}_{\{y_{{\bf i}_{j}} = 0\}}\right)\right) \exp\left( - \frac{a_{n}^{2}}{2} \sum_{j = 0}^{n^{r} - 1} \alpha^{-2j} y_{{\bf i}_{j}}^{2} \right)& \\ &\leq (n^{3/2})^{n^{r}} \exp\left(- \frac{a_{n}^{2}}{2} \inf\limits_{y \in D_{n,k}} \sum_{j = 0}^{n^{r} - 1} \alpha^{-2j} y_{j}^{2}\right) \PP\left(B_{k_{n}}\right),&
\end{align*}
where we used the fact that $1/|y_{{\bf i}_{j}}| \leq n^{3/2}$ if $y_{{\bf i}_{j}} \neq 0$ for the last inequality. From 
\eqref{eq:Ipir0} we obtain
\begin{multline*}
\PP\left(E_{n}({\bf i}) | Y_{n} \geq a_{n},B_{k_{n}}\right) \leq (1/2)^{n-n^{r}} (n^{3/2})^{n^{r}} \exp\left(- \frac{a_{n}^{2}}{2} \inf\limits_{y \in D_{n^{r},k}} \sum_{j = 0}^{n^{r} - 1} \alpha^{-2j} y_{j}^{2}\right) \\ \times \PP\left(B_{k_{n}}\right)\PP\left(Y_{n} \geq a_{n},B_{k_{n}}\right)^{-1}.
\end{multline*}
Finally, \eqref{eq:claim2} follows using \eqref{eq:loyn}.
\end{proof}
As a direct consequence of \eqref{eq:claim1} and \eqref{eq:claim2}, we have, for some positive constant $C$,
\begin{equation}\label{eq:inegTkn}
 \sum_{{\bf i} \in \mathcal{I}^{n}(k)} \PP\left(E_{n}({\bf i}) | Y_{n} \geq a_{n}, B_{k_{n}}\right) \lesssim (2n^{3/2})^{n} e^{- C\frac{a_{n}^{2}}{2}\left(\alpha^{-n^{r}} n^{-3/2}\right)^{2}} + (2n^{3/2})^{n^{r}} (1/2)^{n - n^{r}} R_n. 
\end{equation}
Recalling \eqref{eq:ID}, we obtain from Lemma \ref{lem:opti} that
\begin{equation*}
\inf\limits_{y \in D_{n,k}} \sum_{j = 0}^{n^{r} - 1} \alpha^{-2j} y_{j}^{2} \geq I(\{y  \in \mathbb{R}^{n}; \, \, \,  |y_{k}-v_{k}|\geq \ep\}) > I(H) = 1 - \alpha^{2}.
\end{equation*}
Since $1 - \alpha^{2k_{n} + 2} \rightarrow 1$ when $n \rightarrow \infty,$ we get\begin{equation*}
R_n \leq \exp\left(- C \frac{a_{n}^{2}}{2}\right),
\end{equation*}
 for some positive constant $C$.
From the latter inequality and \eqref{eq:inegTkn}, we get, since $n \rightarrow \infty,$
\begin{equation*}
 \sum_{{\bf i} \in \mathcal{I}^{n}(k)} \PP\left(E_{n}({\bf i}) | Y_{n} \geq a_{n}, B_{k_{n}}\right) \lesssim e^{- C\frac{a_{n}^{2}}{2}\left(\alpha^{-n^{r}} n^{-3/2}\right)^{2}}+ e^{- C \frac{a_{n}^{2}}{2}} \lesssim e^{- C \frac{a_{n}^{2}}{2}}
\end{equation*}
and using \eqref{eq:cvTkn} this ends the proof of \eqref{bornegaussienne}.  
 
\subsection{A law of large numbers}
 \label{ALLN}
 We state here a strong  law of large numbers useful for the proof of local densities in the strictly stable case.
  \begin{prop}  \label{LLNadapted} 
Let $\{\mathcal F_n\}_0^{\infty}$ be a filtration. Let $\{ X_{n,i} : n,i\geq 1\}$ be non-negative real valued
 random variables such that for each $n$,
conditionally on  $\mathcal{F}_n$,
$\{ X_{n,i} : i\geq 1\}$ are independent and identically  r.v. distributed as $X_n$.\\
 Let $\{N_n :  n\geq 1\}$ be non-negative integer valued r.v. such that for each $n$, $N_n$ is $\mathcal F_n$ measurable
and 
 $\liminf_{n\rightarrow \infty} N_{n+1}/N_n>1.$\\
 We assume that  $X_{n,i}$ are   uniformly bounded and $\limsup_{n\rightarrow \infty} \EE(X_{n})/\EE(X_{n+1})< \liminf_{n\rightarrow \infty} N_{n+1}/N_n.$ \\
Then 
$$\frac{1}{N_n\EE(X_n)} \sum_{i=1}^{N_n}  X_{n,i} \stackrel{n\rightarrow \infty}{\longrightarrow} 1 \quad \text{a.s.}.$$
\end{prop}
This result  is a consequence of  a classical law of large numbers where the variables $X_{n,i}$ depend both on $n$ and $i$, see Kurtz \cite{K72} and Athreya Kang \cite{AK98a}.
We prove here that gathering terms in a suitable way allows to deal with the case when $\EE(X_n)\rightarrow 0$.
\begin{proof} We write $u_n= \lfloor 1/\EE(X_n)\rfloor \vee 1$ and
$$ \sum_{i=1}^{N_n}  X_{n,i}=\sum_{k=0}^{V_n-1} Y_{n,k}+R_n,$$
where $V_n=\lfloor N_n/u_n \rfloor$ goes to infinity  and
$$Y_{n,k}=\sum_{i=ku_n}^{(k+1)u_n-1} X_{n,i}, \qquad R_n=\sum_{i=V_nu_n}^{N_n} X_{n,i}.$$
Moreover,   Markov inequality and boundedness of  $X_n$ yields for any $t\geq 0$, $n,k\geq 0$,
$$\PP(Y_{n,k} >t \vert \mathcal F_n)\leq e^{-t}\EE(e^{X_{n}})^{u_n}\leq e^{-t}(1+C\EE(X_n))^{u_n}\leq C' e^{-t}$$
for some constants $C,C'\geq 0$, since $\EE(X_n)$ is bounded.
The right hand side is integrable on $[0,\infty)$ and $\EE(Y_{n,k})$ is bounded and 
$$\liminf_{n\rightarrow \infty} V_{n+1}/V_n\geq  \liminf_{n\rightarrow \infty} N_{n+1}/N_n.  \liminf_{n\rightarrow \infty} \EE(X_{n+1})/\EE(X_n)>1.$$
 We can apply a law of large numbers of  Athreya and Kang (Lemma 1 in \cite{AK98a}) to the family of 
centered and independent variables $\{ Y_{n,k} -\EE(Y_{n,0}) : n\geq 1, k\geq 0\}$ which are stochastically dominated and     get
$$\frac{1}{V_n} \sum_{k=0}^{V_n-1}  (Y_{n,k} -\EE(Y_{n,0})) \stackrel{n\rightarrow \infty}{\longrightarrow} 0 \quad \text{a.s.}$$
and conclude using that $\EE(Y_{n,0})=u_n\EE(X_n)$ and $V_n\sim N_n/u_n$ as $n\rightarrow \infty$ and 
\begin{equation}
\label{partieneg}
 \frac{R_n}{N_n}    =\frac{\sum_{i=V_nu_n}^{N_n} X_{n,i}}{N_n} \stackrel{n\rightarrow \infty}{\longrightarrow} 0 \quad \text{a.s.} 
 \end{equation}
To justify this last limit, one can consider only the  case when $u_n$ goes  geometrically to infinity, since otherwise the boundedness of the r.v. $X_{n,i}$ allows to use a direct domination and Lemma 1 in \cite{AK98a}. When  $u_n$ grows geometrically, using again the law of large numbers in \cite{AK98a} yields
$$\frac{\sum_{i=V_nu_n}^{N_n} X_{n,i}}{u_n} \stackrel{n\rightarrow \infty}{\longrightarrow} 1 \quad \text{a.s.}$$
and the fact that $u_n/N_n\rightarrow 0$ a.s. ensures \eqref{partieneg}  and ends the proof.
\end{proof}
$\newline$
{\bf Acknowledgment.} This work have been supported  by the Chair ``Mod\'elisation Math\'ematique et Biodiversit\'e'' of VEOLIA Environnement-\'Ecole
Polytechnique-MNHN-F.X and ANR ABIM (ANR-16-CE40-0001).

\end{document}